\documentclass[sn-mathphys,Numbered]{sn-jnl}% Math and Physical Sciences Reference Style
\usepackage{graphicx}%
\usepackage{multirow}%
\usepackage{amsmath,amssymb,amsfonts}%
\usepackage{amsthm}%
\usepackage{mathrsfs}%
\usepackage[title]{appendix}%
\usepackage{xcolor}%
\usepackage{textcomp}%
\usepackage{manyfoot}%
\usepackage{booktabs}%
\usepackage{algorithm}%
\usepackage{algorithmicx}%
\usepackage{algpseudocode}%
\usepackage{listings}%
\usepackage{cancel}

\theoremstyle{thmstyleone}%
\theoremstyle{thmstyletwo}%

\theoremstyle{thmstylethree}%

\begin{document}

\title[Article Title]{A Stabilised Semi-Implicit Double-Point Material Point Method for Soil-Water Coupled Problems}


\author*[1]{\fnm{Mian} \sur{Xie}}\email{mian.xie.18@ucl.ac.uk}

\author[2]{\fnm{Pedro} \sur{Navas}}\email{pedro.navas@upm.es}

\author[1]{\fnm{Susana} \sur{L\'opez-Querol}}\email{s.lopez-querol@ucl.ac.uk}

\affil*[1]{\orgdiv{Department of Civil, Environmental and Geomatic Engineering}, \orgname{University College London}, \orgaddress{\street{Gower Street}, \city{London}, \postcode{WC1E 6BT}, \country{UK}}}

\affil[2]{\orgdiv{ETSI Caminos Canales y Puertos}, \orgname{Universidad Polit\'ecnica de Madrid}, \orgaddress{\street{Calle Profesor Aranguren 3}, \city{Madrid}, \postcode{28040}, \country{Spain}}}


\abstract{This paper presents a novel semi-implicit two-phase double-point Material Point Method (MPM) for modelling large deformation in geotechnical engineering problems. To overcome the computational limitations of explicit methods, we develop a semi-implicit approach that eliminates time step dependency in soil-water coupled problems. Unlike existing single-point methods that use one set of material points to represent soil-water mixtures, our approach employs two distinct sets of material points to model soil and water phases separately. To address MPM's inherent stress oscillations, we introduce a stabilisation technique based on the modified F-bar method.
Through validation against Terzaghi's one-dimensional consolidation theory, one-dimensional large deformation consolidation, and large deformation slope stability studies, our method demonstrates superior performance. Further testing with the hyperelastic Nor-Sand constitutive model in landslide simulations reveals that the double-point approach produces significantly more reliable results than single-point methods, particularly for dilatant soils. Notably, while implementing two sets of material points, our method incurs only a 15\% increase in computational cost while achieving markedly improved accuracy. These findings establish the double-point MPM as a robust and efficient approach for analysing large deformation geotechnical problems under fully saturated conditions.
}

\keywords{Material Point Method, Incremental fractional-step method, Finite strain, Stabilisation, Slope stability}



\maketitle

\hspace{1cm}

\section{Introduction}\label{sec1}
The material point method (MPM), originally developed by \citet{Sulsky1994}, has become a popular numerical tool to solve large deformation geotechnical problems in recent years \citep{Soga2016,Augarde2021}. Frequently, coupled analysis in saturated soils is of crucial importance in geotechnical engineering, but adds more complexity to the analysis than only dealing with a single phase.

Existing soil-water coupled (two-phase) MPMs can be classified into two categories: the single-point and the double-point approaches. In a two-phase, single-point MPM \citep{Zhang2009,Alonso2010,Zabala2011,Jassim2013,Iaconeta2019,Zheng2022,Kularathna2021,Yuan2023}, one set of material points is used to model the soil-water mixture. The same set of material points carries information for both soil and water phases, for example, the stress (soil phase) and the pore pressure (water phase). In this case, the properties of the water phase, such as pore pressure, are attached to the soil material points, and the water phase follows the movement of the soil skeleton \citep{Yerro2015}. Therefore, the relative acceleration between the soil and water phases is normally ignored for simplicity \citep{Soga2016}. 

In a two-phase, double-point MPM \citep{Abe2014,Bandara2015,Liu2017,Yamaguchi2020}, two sets of material points are used to model soil and water independently. The soil and water material points carry their own information. In addition, the movement of the water material points is no longer attached to the soil skeleton. The superposition of the soil and water material points represents the soil-water mixture. The stored information is mapped from the soil and water material points to the same grid nodes where the governing equations of the mixture are resolved. Compared with the single-point approach, the obvious disadvantage is that it requires a higher computational effort because of the additional set of material points. The advantage, however, is that the double-point approach is more suitable to model problems such as seepage, erosion, and submarine landslides, in which the movement of the soil and fluid phases is considerably different. Furthermore, \citet{Soga2016} claim that it is important to use a double-point approach when the relative acceleration between the soil and water is high. However, they neither provide a clear recommendation for this limiting relative acceleration, nor a numerical example to compare the differences of both procedures.

%In the two-phase MPMs mentioned above, the air phase is ignored. Therefore, they are not capable of modelling unsaturated soils the pores of which are filled with both air and water, and the suction plays an important role. However, it is possible to modify a two-phase formulation to accommodate suction \citep{Bandara2016,Wang2018,Martinelli2021,Feng2021}. In these approaches, the gas pressure and density are ignored. Alternatively, the unsaturated soil can be modelled with a three-phase MPM (that is, a soil-water-air coupled MPM) \citep{Yerro2015,Zheng2022b}. All existing three-phase MPMs \citep{Yerro2015,Zheng2022b} use a single-point approach, due to computational efficiency reasons. Additionally, the soil's constitutive model has to be modified as an unsaturated version.   

\citet{Zheng2022} reviewed the existing MPM studies for large deformation problems in fluid-saturated porous media. From \citet{Zheng2022}, we can see that only very basic constitutive models, such as elastic, Drucker-Prager, and Mohr-Coulomb, are used in these {soil-water coupled} analyses. Exceptions are the works of \citet{Ceccato2016} and \citet{Martinelli2022}, where the modified Cam-Clay model is employed to study penetration problems. Recently, \citet{Carluccio2023} and \citet{Lino2023}, from the Anura3D community, studied the problem of landslides using the Ta-Ger and SANISAND models, respectively. Therefore, the performance of advanced constitutive models under large deformation {soil-water coupled problems} remains to be studied in depth. {Recent studies of the pile installation in the single-phase MPM also emphasise the importance of using an advanced constitutive model to study a large deformation problem \citep{Martinelli2021b,Martinelli2022b,Fetrati2024,Galavi2024}.}

From all these existing {soil-water coupled} MPM studies in the literature that use advanced soil models \citep{Ceccato2016,Martinelli2022,Carluccio2023,Lino2023}, we can conclude that an explicit dynamic {coupled} MPM has been implemented in all cases. An advanced soil model may not be feasible in an implicit {coupled} MPM formulation at the moment due to convergence issues under large deformation ranges, because implicit algorithms require a stiffness matrix derived from the constitutive model and the derivation of this matrix may be challenging for those models under large deformation. Further research is yet needed in this line. 

By contrast, it is well known that the stiffness matrix is not needed in an explicit MPM. However, a fully explicit soil-water coupled MPM requires a very small time step compared with an implicit approach. The critical time step of an explicit approach can be determined by the well established Courant–Friedrichs–Lewy (CFL) condition \citep{Courant1967}, according to which the critical time step is inversely proportional to the bulk moduli of the mixture. Therefore, in this case, the water phase dominates the critical time step because its bulk modulus is much higher than that of the soil, resulting in a very small time step. Furthermore, \citet{Mieremet2016} found that the critical time step is also related to the permeability for an explicit two-phase method. Under low-permeability conditions, the actual critical time step is even smaller than the one given by the CFL condition.

\citet{Yamaguchi2020} developed a semi-implicit, double-point MPM, based on the fractional-step method proposed by \citet{Kularathna2017}. Semi-implicit means that the soil phase is explicit, while the water phase is solved implicitly. Therefore, the critical time step of a semi-implicit MPM is dominated by the CFL condition of the soil phase, resulting in a significant improvement in computational efficiency compared with an explicit, two-phase MPM. In addition, the pore pressure is more stable in the fractional-step method because it satisfies the inf-sup condition. However, \citet{Kularathna2021} realised that the dissipation of the pore pressure in this fractional-step method depends on the time step, resulting in poor performance in geotechnical applications. {For the same reason, the critical time step in the approach \citet{Yamaguchi2020} proposed requires a very small time step (far smaller than the CFL condition) in the case of low permeability to maintain numerical stability.} Therefore, \citet{Kularathna2021} developed a semi-implicit single-point MPM based on an incremental fractional-step method to remove the dependency of the time step. \citet{Yuan2023} significantly improved the computational efficiency of the work of \citet{Kularathna2021} by introducing a novel representation of the drag force between soil and water phases. However, the work of \citet{Yuan2023} is still based on a single-point formulation. The model developed by \citet{Pan2021} is also a semi-implicit approach, but the water phase is solved implicitly using the Eulerian finite element method (FEM). This MPM-FEM coupled approach has advantages in the definition of the free water surface and the treatment of inflow and outflow boundaries \citep{Pan2021}. Therefore, it is more appropriate for problems such as submarine landslides in which free water dominates the problem. Nevertheless, a very basic constitutive model, either Mohr–Coulomb or Drucker–Prager, is used in these studies. The performances of the MPMs have to be investigated using a more advanced soil constitutive model.               

In the current study, we further derive a semi-implicit, two-phase, double-point MPM from the single-point method proposed by \citet{Yuan2023}. The novelty of this research compared with \citet{Yuan2023} is that the weak forms are discretised using two sets of material points to represent the soil and water phases separately. {Also, the higher-order B-spline shape function is used in this research.} The soil and water phases of the proposed method are stabilised with a new method that is developed to fully eliminate volumetric locking, based on the modified F-bar method proposed by \citet{Xie2023b}. {The newly developed modified F-bar method is specially tailored for higher-order shape function and advanced constitutive model.} Nor-Sand constitutive model is implemented in the proposed method to highlight the different results obtained with the single-point and double-point approaches, and the importance of using a double-point MPM is concluded. 

The rest of this paper is organised as follows: Section \ref{sec2} introduces some basic concepts, assumptions, and governing equations for the proposed semi-implicit two-phase double-point MPM. Then, the time discretisation of these governing equations, based on the incremental fractional-step method, is presented in Section \ref{sec3}. In Section \ref{sec4}, the weak form, spatial discretisation, and numerical algorithm for the newly developed semi-implicit, two-phase, double-point MPM are presented. The finite strain constitutive models and stabilisation methods are also discussed in Section \ref{sec4}. Section \ref{sec5} presents some numerical models for validation and study purposes, including Terzaghi’s one-dimensional consolidation, {one-dimensional large deformation consolidation,} slope stability with the Mohr-Coulomb model, as well as landslides with the Nor-Sand model. Finally, Section \ref{sec6} concludes the important findings and limitations of this research.

\section{Governing equations}\label{sec2}

\subsection{Some basics}
In the soil-water mixture employed hereinafter, it is assumed that the soil is fully saturated. As previously mentioned, we use two sets of material points to represent the soil and water phases separately, and the mixture is considered to be the superposition of the soil and water material points. Both soil and water phases follow the Lagrangian descriptions. In this paper, the subscripts $s$ and $w$ represent the soil and water phases, respectively. Therefore, the coordinates of the soil and water material points in the reference configuration (undeformed) are written as $\pmb{X}_{s}$ and $\pmb{X}_{w}$, respectively. Similarly, the spatial coordinates (in current or deformed configuration) for soil and water material points are $\pmb{x}_{s}$ and $\pmb{x}_w$, respectively. The reference and spatial coordinates of a material point are related by a motion $\chi_{\alpha}$:
\begin{equation}
    \pmb{x}_{\alpha} = \chi_{\alpha}(\pmb{X}_{\alpha},\,t)\,,
\end{equation}
where the subscript $\alpha = s,\,w$ represents either soil or water, and $t$ denotes the time elapsed from the reference to the current configuration. For the updated Lagrangian formulation adopted in this research, the last converged configuration at time $t$ is taken as the reference configuration for the next time step $t+1$.

The velocity $\pmb{v}_{\alpha}$ and acceleration $\pmb{a}_{\alpha}$ for soil or water 
material points are defined as
\begin{equation}
    \pmb{v}_{\alpha} = \frac{\text{d}_\alpha\pmb{x}_{\alpha}}{\text{d}t} \,;
\end{equation}

\begin{equation}
    \pmb{a}_{\alpha} = \frac{\text{d}_\alpha\pmb{a}_{\alpha}}{\text{d}t} \,.
\end{equation}

According to the mixture theory, the density of mixture $\rho$ is defined as
\begin{equation}
    \rho = (1-n)\,\rho_s + n\,\rho_w \,,
\end{equation}
where $n$ is the porosity of soil, $\rho_s$ is the density of soil grains, and $\rho_w$ is the density of the water. For simplicity, both soil grains and water densities remain constant during the analysis, while porosity is updated throughout the simulation. 

In this research, Terzaghi's effective stress principle is adopted. Therefore, the total stress of soil $\pmb{\sigma}_s$ is defined as
\begin{equation}
	\pmb{\sigma}_s = \pmb{\sigma}' + (1-n)\,p\,\pmb{I} \,,
\end{equation}
where $\pmb{\sigma}'$ is the effective stress of the soil skeleton, $p$ is the pore water pressure, and $\pmb{I}$ is the identity tensor. The total stress of water $\pmb{\sigma}_w$ is defined as
\begin{equation}
	\pmb{\sigma}_w = n\,p\,\pmb{I} \,.
\end{equation}
In this paper, to ease the numerical implementation, tension is always considered with a positive sign convention for both soil stress and pore pressure.

\subsection{Conservation of mass}
In this research, both the soil grains and the water are assumed as incompressible, and it can be mathematically written as 
\begin{equation}
    \frac{\text{d}_{\alpha}\rho_{\alpha}}{\text{d}t} = 0 \,.
\end{equation}
By taking into account the assumption of incompressibility, the mass balance equations for the soil and water phases can be written as
\begin{equation}\label{eq:mass_balance_soil}
    \frac{\text{d}_s(1-n)}{\text{d}t} + (1-n)\,\nabla\cdot\pmb{v}_s = 0 \,,
\end{equation}
and
\begin{equation}\label{eq:mass_balance_water}
    \frac{\text{d}_w n}{\text{d}t} + n\,\nabla\cdot\pmb{v}_w = 0 \,,
\end{equation}
where
\begin{equation}
    \frac{\text{d}_{\alpha} (\Box)}{\text{d}t} = \frac{\partial(\Box)}{\partial t} + \pmb{v}_{\alpha}\cdot\nabla(\Box)\,.
\end{equation}
By summing up Eqs.~(\ref{eq:mass_balance_soil}) and~(\ref{eq:mass_balance_water}), the mass balance equation for the mixture can be obtained as
\begin{equation}\label{eq:mass_balance_mixture}
    \nabla\cdot[(1-n)\,\pmb{v}_s + n\,\pmb{v}_w] = 0 \,.
\end{equation}
We assume that the porosity gradient is very small in this research. Therefore, Eq.~(\ref{eq:mass_balance_mixture}) can be simplified as
\begin{equation}\label{eq:mass_balance_mixture_simplified}
    (1-n)\,\nabla\cdot\pmb{v}_s + n\,\nabla\cdot\pmb{v}_w = 0 \,.
\end{equation}

\subsection{Conservation of momentum}
The momentum balance equations for the soil and water phases can be written as
\begin{equation}\label{eq:momentum_balance_soil}
    (1-n)\,\rho_s\,\pmb{a}_s = \nabla\cdot[\pmb{\sigma}' + (1-n)\,p\,\pmb{I}] + (1-n)\,\rho_s\,\pmb{b} + \hat{\pmb{p}}  \,,
\end{equation}
and
\begin{equation}\label{eq:momentum_balance_water}
    n\,\rho_w\,\pmb{a}_w = \nabla\cdot(n\,p\,\pmb{I}) + n\,\rho_w\,\pmb{b} - \hat{\pmb{p}} \,,
\end{equation}
where $\pmb{b}$ is the body force and $\hat{\pmb{p}}$ represents the momentum interaction between soil skeleton and water. By assuming a laminar flow (i.e. linear Darcy’s law) and neglecting the buoyancy term (i.e. $p\,\nabla\,n$), the momentum interaction $\hat{\pmb{p}}$ can be written as
\begin{equation}
    \hat{\pmb{p}} = \frac{n^2\,\rho_w\,g}{k}\cdot(\pmb{v}_w - \pmb{v}_s) \,,
\end{equation}
where $g\approx9.81\,m/s^2$ is the gravity acceleration and $k$ is the hydraulic conductivity with units of $[m/s]$.

By adding Eqs.~(\ref{eq:momentum_balance_soil}) and~(\ref{eq:momentum_balance_water}), the momentum balance equation for the mixture can be obtained as
%\begin{multline}\label{eq:momentum_balance_mixture}
%    (1-n)\,\rho_s\,\pmb{a}_s + n\,\rho_w\,\pmb{a}_w = \\
%    \nabla\cdot(\pmb{\sigma}' + p\,\pmb{I}) + (1-n)\,\rho_s\,\pmb{b} + n\,\rho_w\,\pmb{b} \,.
%\end{multline}
\begin{equation}\label{eq:momentum_balance_mixture}
    (1-n)\,\rho_s\,\pmb{a}_s + n\,\rho_w\,\pmb{a}_w = \nabla\cdot(\pmb{\sigma}' + p\,\pmb{I}) + (1-n)\,\rho_s\,\pmb{b} + n\,\rho_w\,\pmb{b} \,.
\end{equation}

\section{Time discretisation and fractional-step method}\label{sec3}
Equations~(\ref{eq:momentum_balance_mixture}),~(\ref{eq:momentum_balance_water}) and~(\ref{eq:mass_balance_mixture_simplified}) are the governing equations used in this research. Their time-discretised forms can be written as
%\begin{multline}\label{eq:time_discretised_mixture}
%    (1-n)\,\rho_s\,\pmb{a}^{t+1}_s + n\,\rho_w\,\pmb{a}^{t+1}_w = \\
%    \nabla\cdot{\pmb{\sigma}'}^t + \nabla\cdot(p^{t+1}\pmb{I}) + (1-n)\,\rho_s\,\pmb{b} + n\,\rho_w\,\pmb{b} \,,
%\end{multline}
\begin{equation}\label{eq:time_discretised_mixture}
    (1-n)\,\rho_s\,\pmb{a}^{t+1}_s + n\,\rho_w\,\pmb{a}^{t+1}_w = \nabla\cdot{\pmb{\sigma}'}^t + \nabla\cdot(p^{t+1}\pmb{I}) + (1-n)\,\rho_s\,\pmb{b} + n\,\rho_w\,\pmb{b} \,,
\end{equation}
\begin{equation}\label{eq:time_discretised_water}
    n\,\rho_w\,\pmb{a}^{t+1}_w = n\,\nabla\cdot(p^{t+1}\pmb{I}) + n\,\rho_w\,\pmb{b} - \hat{\pmb{p}} \,,
\end{equation}
%\vskip 1mm
\begin{equation}\label{eq:time_discretised_mass}
    (1-n)\,\nabla\cdot\pmb{v}^{t+1}_s + n\,\nabla\cdot\pmb{v}^{t+1}_w = 0 \,,
\end{equation}
where the superscript $t+1$ denotes the corresponding variable at the next configuration. 

The idea of the fractional-step method is to introduce an intermediate velocity $\pmb{v}^{*}_{\alpha}$ (i.e. a velocity field between $\pmb{v}^{t}_{\alpha}$ and $\pmb{v}^{t+1}_{\alpha}$) to split the governing equations. As a result, the calculation of pore pressure and velocity can be decoupled. In general, the fractional step method is a prediction-correction scheme. First, the intermediate velocity $\pmb{v}^{*}_{\alpha}$ is predicted based on the variables in the current configuration. Then, the pore pressure for the next configuration $p^{t+1}$ can be implicitly calculated based on the intermediate velocity $\pmb{v}^{*}_{\alpha}$. Finally, the true velocity for the next configuration $\pmb{v}^{t+1}_{\alpha}$ is updated based on pressure $p^{t+1}$ and intermediate velocity $\pmb{v}^{*}_{\alpha}$. This intermediate velocity field does not necessarily satisfy the incompressibility constraint. Therefore, the pore pressure instability due to incompressibility can be solved by introducing this intermediate velocity field. 

By introducing the intermediate velocity, the acceleration on the next configuration $\pmb{a}^{t+1}_{\alpha}$ can be split as
\begin{equation}\label{eq:fractional_step_acc}
    \pmb{a}^{t+1}_{\alpha} = \frac{\pmb{v}^{t+1}_{\alpha}-\pmb{v}^{t}_{\alpha}}{\Delta t} = \frac{\pmb{v}^{t+1}_{\alpha}-\pmb{v}^{*}_{\alpha}}{\Delta t} + \frac{\pmb{v}^{*}_{\alpha}-\pmb{v}^{t}_{\alpha}}{\Delta t} \,.
\end{equation}

By introducing Eq. (\ref{eq:fractional_step_acc}) into Eq. (\ref{eq:time_discretised_mixture}), this time-discretised momentum balanced equation for the mixture can be split into two parts:
%\begin{multline}\label{eq:soil_phase_part1}
%    (1-n)\,\rho_s\frac{\pmb{v}^{*}_s-\pmb{v}^{t}_s}{\Delta t} + n\,\rho_w\frac{\pmb{v}^{*}_w-\pmb{v}^{t}_w}{\Delta t} = \\
%    \nabla\cdot{\pmb{\sigma}'}^t + \nabla\cdot(p^{t}\pmb{I}) + (1-n)\,\rho_s\,\pmb{b} + n\,\rho_w\,\pmb{b} \,,
%\end{multline}
\begin{equation}\label{eq:soil_phase_part1}
    (1-n)\,\rho_s\frac{\pmb{v}^{*}_s-\pmb{v}^{t}_s}{\Delta t} + n\,\rho_w\frac{\pmb{v}^{*}_w-\pmb{v}^{t}_w}{\Delta t} = \nabla\cdot{\pmb{\sigma}'}^t + \nabla\cdot(p^{t}\pmb{I}) + (1-n)\,\rho_s\,\pmb{b} + n\,\rho_w\,\pmb{b} \,,
\end{equation}
and
%\begin{multline}\label{eq:soil_phase_part2}
%    (1-n)\,\rho_s\frac{\pmb{v}^{t+1}_s-\pmb{v}^{*}_s}{\Delta t} + n\,\rho_w\frac{\pmb{v}^{t+1}_w-\pmb{v}^{*}_w}{\Delta t} = \\
%    \nabla(p^{t+1} - p^t)  \,.
%\end{multline}
\begin{equation}\label{eq:soil_phase_part2}
    (1-n)\,\rho_s\frac{\pmb{v}^{t+1}_s-\pmb{v}^{*}_s}{\Delta t} + n\,\rho_w\frac{\pmb{v}^{t+1}_w-\pmb{v}^{*}_w}{\Delta t} = \nabla(p^{t+1} - p^t)  \,.
\end{equation}

Similarly, the time-discretised momentum balanced equation for the water phase (i.e. Eq. (\ref{eq:time_discretised_water})) can also be split into two parts:
\begin{equation}\label{eq:water_phase_part1_org}
    n\,\rho_w\frac{\pmb{v}^{*}_w-\pmb{v}^{t}_w}{\Delta t} = n\,\nabla\cdot(p^{t}\pmb{I}) + n\,\rho_w\,\pmb{b} - \hat{\pmb{p}} \,,
\end{equation}
and
\begin{equation}\label{eq:water_phase_part2}
    n\,\rho_w\frac{\pmb{v}^{t+1}_w-\pmb{v}^{*}_w}{\Delta t} = n\nabla(p^{t+1} - p^t)  \,.
\end{equation}
In Eq. (\ref{eq:water_phase_part1}), we adopt the momentum interaction term proposed by \citet{Yuan2023}, which has a form of 
\begin{equation}\label{eq:momentum_term}
    \hat{\pmb{p}} = \frac{n^2\,\rho_w\,g}{k}\cdot(\pmb{v}^{*}_w - \pmb{v}^t_s) \,.
\end{equation}
Substituting Eq. (\ref{eq:momentum_term}) into Eq. (\ref{eq:water_phase_part1_org}) and rearranging the equation, we can have
%\begin{multline}\label{eq:water_phase_part1}
%    n\,\rho_w\pmb{v}^*_w(\frac{1}{\Delta t}+\frac{n\,g}{k}) =\\
%    n\,\nabla p^t + n\,\rho_w\pmb{b} + \frac{n^2\rho_w g}{k}\pmb{v}^t_s + \frac{n\,\rho_w}{\Delta t}\pmb{v}^t_w
%\end{multline}
\begin{equation}\label{eq:water_phase_part1}
    n\,\rho_w\pmb{v}^*_w\left(\frac{1}{\Delta t}+\frac{n\,g}{k}\right) = n\,\nabla p^t + n\,\rho_w\pmb{b} + \frac{n^2\rho_w g}{k}\pmb{v}^t_s + \frac{n\,\rho_w}{\Delta t}\pmb{v}^t_w \,.
\end{equation}
The advantage of using Eq. (\ref{eq:momentum_term}) is that the intermediate velocity of water $\pmb{v}^{*}_w$ is the only unknown variable in Eq. (\ref{eq:water_phase_part1}). As a result, the intermediate velocity of water $\pmb{v}^{*}_w$ can be easily obtained using an explicit approach. With $\pmb{v}^{*}_w$ in hand, the intermediate velocity of soil $\pmb{v}^{*}_s$ becomes the only unknown variable in Eq. (\ref{eq:soil_phase_part1}), which can also be solved explicitly. In the original method proposed by \citet{Kularathna2021}, the intermediate velocities of soil and water must be solved simultaneously, resulting in significant computational effort. 

Rearranging Eqs. (\ref{eq:soil_phase_part2}) and (\ref{eq:water_phase_part2}), the velocities in the next time step, $\pmb{v}^{t+1}_\alpha$, are obtained. Substituting these velocities $\pmb{v}^{t+1}_\alpha$ into the mass balance equation for the mixture (Eq. (\ref{eq:time_discretised_mass})), we can obtain the Poisson equation
%\begin{multline}\label{eq:possion_eqn}
%    (1-n)\,\nabla\cdot\pmb{v}^{*}_s + n\,\nabla\cdot\pmb{v}^{*}_w \\
%    + \Delta t(\frac{1-n}{\rho_s}+\frac{n}{\rho_w})\,\nabla^2(p^{t+1}-p^t) = 0 \,.
%\end{multline}
\begin{equation}\label{eq:possion_eqn}
    (1-n)\,\nabla\cdot\pmb{v}^{*}_s + n\,\nabla\cdot\pmb{v}^{*}_w   + \Delta t\left(\frac{1-n}{\rho_s}+\frac{n}{\rho_w}\right)\,\nabla^2(p^{t+1}-p^t) = 0 \,.
\end{equation}
As we can see in Eq.~(\ref{eq:possion_eqn}), the only unknown is the pore pressure increment $(p^{t+1}-p^t)$, which can be obtained implicitly.

Finally, the velocity of water and soil phases can be updated using Eqs.~(\ref{eq:water_phase_part2}) and~(\ref{eq:soil_phase_part2}), respectively.

\section{Coupled MPM formulation and algorithm}\label{sec4}

\subsection{Weak form of the governing equations}\label{subsec:weak_form}
The weak form of an MPM formulation is the same as the FEM, which is obtained using the Galerkin method. In this method, a test function is multiplied on both sides of the governing equations, which are then integrated over the current configuration $\Omega_{\alpha}$ of the soil or water phase. Also, integration by parts using the Green-Gauss divergence theorem is applied to the stress terms. Finally, the weak form of governing equations for the mixture (Eqs.~(\ref{eq:soil_phase_part1}) and (\ref{eq:soil_phase_part2})) are written as
%\begin{multline}
%    \int_{\Omega_{s}}^{} \delta \pmb{v}_s \cdot(1-n)\rho_s \frac{\pmb{v}^*_s-\pmb{v}^t_s}{\Delta t} \, dV \\ + 
%    \int_{\Omega_{w}}^{} \delta \pmb{v}_w \cdot n\rho_w \frac{\pmb{v}^*_w-\pmb{v}^t_w}{\Delta t} \, dV = \\
%    -\int_{\Omega_{s}}^{} \nabla\delta \pmb{v}_s : {\pmb{\sigma}'}^t \, dV + \int_{\partial\Omega_{\tau}}^{} \delta \pmb{v}_s \cdot \Bar{\pmb{\tau}} \, dS \\ -\int_{\Omega_{w}}^{} \nabla\delta \pmb{v}_w : (p^t\pmb{I}) \, dV \\
%    + \int_{\Omega_{s}}^{} \delta \pmb{v}_s \cdot (1-n)\rho_s \, \pmb{b} \, dV + \int_{\Omega_{w}}^{} \delta \pmb{v}_w \cdot n\rho_w \, \pmb{b} \, dV  \,, 
%\end{multline}
\begin{multline}
    \int_{\Omega_{s}}^{} \delta \pmb{v}_s \cdot(1-n)\rho_s \frac{\pmb{v}^*_s-\pmb{v}^t_s}{\Delta t} \, dV  + 
    \int_{\Omega_{w}}^{} \delta \pmb{v}_w \cdot n\rho_w \frac{\pmb{v}^*_w-\pmb{v}^t_w}{\Delta t} \, dV = \\ -\int_{\Omega_{s}}^{} \nabla\delta \pmb{v}_s : {\pmb{\sigma}'}^t \, dV + \int_{\partial\Omega_{\tau}}^{} \delta \pmb{v}_s \cdot \Bar{\pmb{\tau}} \, dS \\
    -\int_{\Omega_{w}}^{} \nabla\delta \pmb{v}_w : (p^t\pmb{I}) \, dV 
    + \int_{\Omega_{s}}^{} \delta \pmb{v}_s \cdot (1-n)\rho_s \, \pmb{b} \, dV + \int_{\Omega_{w}}^{} \delta \pmb{v}_w \cdot n\rho_w \, \pmb{b} \, dV  \,, 
\end{multline}
and
%\begin{multline}
%    \int_{\Omega_{s}}^{} \delta \pmb{v}_s \cdot(1-n)\rho_s \frac{\pmb{v}^{t+1}_s-\pmb{v}^*_s}{\Delta t} \, dV \\ + 
%    \int_{\Omega_{w}}^{} \delta \pmb{v}_w \cdot n\rho_w \frac{\pmb{v}^{t+1}_w-\pmb{v}^*_w}{\Delta t} \, dV = \\
%    \int_{\Omega_{s}}^{} (1-n)\nabla\delta \pmb{v}_s : (p^{t+1} - p^t) \, dV \\
%    +\int_{\Omega_{w}}^{} n\nabla\delta \pmb{v}_w : (p^{t+1} - p^t) \, dV \,,
%\end{multline}
\begin{multline}
    \int_{\Omega_{s}}^{} \delta \pmb{v}_s \cdot(1-n)\rho_s \frac{\pmb{v}^{t+1}_s-\pmb{v}^*_s}{\Delta t} \, dV  + 
    \int_{\Omega_{w}}^{} \delta \pmb{v}_w \cdot n\rho_w \frac{\pmb{v}^{t+1}_w-\pmb{v}^*_w}{\Delta t} \, dV = \\
    \int_{\Omega_{s}}^{} (1-n)\nabla\delta \pmb{v}_s : (p^{t+1} - p^t) \, dV 
    +\int_{\Omega_{w}}^{} n\nabla\delta \pmb{v}_w : (p^{t+1} - p^t) \, dV \,,
\end{multline}
where $\Bar{\pmb{\tau}} = \pmb{\sigma}' \cdot \pmb{n}$ is the prescribed traction acting on the soil phase's boundary $\partial\Omega_{\tau}$; $\pmb{n}$ is the unit vector normal to the boundary where the force is prescribed; $dS$ is the surface element where the force is prescribed; $\delta\pmb{v}_{\alpha}$ is the test function. The choice of test function is arbitrary, but the test function becomes zero at the nodes where displacements are prescribed.

Following the same procedure, the weak forms of the governing equation for the water phase (Eqs.~(\ref{eq:water_phase_part1}) and~(\ref{eq:water_phase_part2})) are written as
%\begin{multline}
%    \int_{\Omega_{w}}^{} \delta \pmb{v}_w \cdot n\,\rho_w\pmb{v}^*_w(\frac{1}{\Delta t}+\frac{n\,g}{k}) \, dV =\\
%    - \int_{\Omega_{w}}^{} n\, \nabla\delta \pmb{v}_w : (p^t\pmb{I}) \, dV  + \int_{\Omega_{w}}^{} \delta \pmb{v}_w \cdot n\,\rho_w\pmb{b} \, dV \\
%    + \int_{\Omega_{w}}^{} \delta \pmb{v}_w \cdot \frac{n^2\rho_w g}{k}\pmb{v}^t_s \, dV + \int_{\Omega_{w}}^{} \delta \pmb{v}_w \cdot \frac{n\,\rho_w}{\Delta t}\pmb{v}^t_w \, dV \,,
%\end{multline}
\begin{multline}
    \int_{\Omega_{w}}^{} \delta \pmb{v}_w \cdot n\,\rho_w\pmb{v}^*_w\left(\frac{1}{\Delta t}+\frac{n\,g}{k}\right) \, dV =
    - \int_{\Omega_{w}}^{} n\, \nabla\delta \pmb{v}_w : (p^t\pmb{I}) \, dV  + \int_{\Omega_{w}}^{} \delta \pmb{v}_w \cdot n\,\rho_w\pmb{b} \, dV \\
    + \int_{\Omega_{w}}^{} \delta \pmb{v}_w \cdot \frac{n^2\rho_w g}{k}\pmb{v}^t_s \, dV + \int_{\Omega_{w}}^{} \delta \pmb{v}_w \cdot \frac{n\,\rho_w}{\Delta t}\pmb{v}^t_w \, dV \,,
\end{multline}
and
%\begin{multline}
%    \int_{\Omega_{w}}^{} \delta \pmb{v}_w \cdot n\,\rho_w\frac{\pmb{v}^{t+1}_w-\pmb{v}^{*}_w}{\Delta t}\,dV =\\
%    \int_{\Omega_{w}}^{} n\,\nabla\delta \pmb{v}_w : (p^{t+1} - p^t)\,dV  \,.
%\end{multline}
\begin{equation}
    \int_{\Omega_{w}}^{} \delta \pmb{v}_w \cdot n\,\rho_w\frac{\pmb{v}^{t+1}_w-\pmb{v}^{*}_w}{\Delta t}\,dV =
    \int_{\Omega_{w}}^{} n\,\nabla\delta \pmb{v}_w : (p^{t+1} - p^t)\,dV  \,.
\end{equation}

Finally, the weak form of the Poisson equation (Eq.~(\ref{eq:possion_eqn})) is written as
%\begin{multline}
%    \int_{\Omega_{s}}^{} \delta p\,(1-n)\,\nabla\cdot\pmb{v}^{*}_s\,dV + \int_{\Omega_{w}}^{} \delta p\,n\,\nabla\cdot\pmb{v}^{*}_w \,dV \\
%    - \int_{\Omega_{s}}^{}\Delta t \nabla\delta p\frac{1-n}{\rho_s}\nabla(p^{t+1}-p^t)\,dV \\
%    - \int_{\Omega_{w}}^{}\Delta t \nabla\delta p\frac{n}{\rho_w}\nabla(p^{t+1}-p^t)\,dV = 0 \,,
%\end{multline}
\begin{multline}\label{eq:possion_eqn_weak}
    \int_{\Omega_{s}}^{} \delta p\,(1-n)\,\nabla\cdot\pmb{v}^{*}_s\,dV + \int_{\Omega_{w}}^{} \delta p\,n\,\nabla\cdot\pmb{v}^{*}_w \,dV 
    - \int_{\Omega_{s}}^{}\Delta t \nabla\delta p\frac{1-n}{\rho_s}\nabla(p^{t+1}-p^t)\,dV \\
    - \int_{\Omega_{w}}^{}\Delta t \nabla\delta p\frac{n}{\rho_w}\nabla(p^{t+1}-p^t)\,dV = 0 \,,
\end{multline}
where $\delta p$ is the test function for the Poisson equation.

\subsection{Spatial discretisation}

In this research, we use two sets of material points to discretise the soil and water continuum bodies. In the initial configuration $\Omega^0_{\alpha}$, each material point occupies a volume $V^0_{\alpha p}$ in space, which represents a portion of the corresponding continuum body. The subscript $p$ represents the variable associated with a material point. The mass of soil in the initial configuration can be obtained by summing up all the contributions of soil material points:
\begin{equation}\label{eq:msp}
    m^0_{s} = \int_{\Omega^0_{s}}^{} (1-n^0)\rho_s \, dV \approx \sum_{s p=1}^{{N}_{s p}}  (1-n^0)\rho_s \, V^0_{sp} = \sum_{s p=1}^{{N}_{s p}}  m^0_{s p} \,,
\end{equation}
where $n^0$ is the initial porosity; $\Omega^0_{s}$ is the initial domain of the soil phase; ${N}_{s p}$ is the total number of soil material points; and $m^0_{s p}$ is the initial mass of soil material point.
In this study, the total masses of soil and water are conserved during the calculation:
\begin{equation}\label{eq:conv_mass}
    m_{\alpha} = m^{t}_{\alpha} \equiv m^0_{\alpha} \,.
\end{equation} 
Therefore, the mass of a {soil} material point remains constant in time:
\begin{equation}
    m_{s p} = m^{t}_{s p} \equiv m^0_{s p} \,.
\end{equation}
However, the volume of each soil material point is updated during the calculation. 

{During the simulation, the porosity $n$ is updated, changing the effective masses of soil and water phases. In a double-point formulation, a water material point is allowed to travel through the porous media. Therefore, the volume of a material point is updated based on the changing of porosity:
\begin{equation}
    V_{wp} = n^{0}/n \, V_{wp}^{0} \,.
\end{equation}
The mass of water can be obtained by 
\begin{equation}\label{eq:mwp}
    m_{w} = \int_{\Omega_{w}}^{} n\,\rho_w \, dV \approx \sum_{w p=1}^{{N}_{w p}} n\,\rho_w \, V_{wp} =  \sum_{w p=1}^{{N}_{w p}} n^0\,\rho_w \, V_{wp}^0 \equiv m^0_{w} \,,
\end{equation}
where $\Omega_{w}$ is the domain of the water phase, and ${N}_{w p}$ is the total number of water material points.}

MPM has very similar spatial discretisation as FEM. In MPM, the variables at a material point location can be interpolated from the corresponding nodal values. We can see that the material point is equivalent to the quadrature point in FEM. As in FEM, the shape function (or basis function) is used to interpolate the velocity and pressure fields. Therefore, the velocity, intermediate velocity, and pressure at the location of the material point, can be approximated as
\begin{equation}\label{eq:vI2vp}
    \pmb{v}^t_{\alpha p}(\pmb{x}_{\alpha}) = \sum_{I=1}^{{N}_{n}} {S}_{I\alpha p}\,\pmb{v}^t_{\alpha I}\,,
\end{equation}
\begin{equation}
    \pmb{v}^*_{\alpha p}(\pmb{x}_{\alpha}) = \sum_{I=1}^{{N}_{n}} {S}_{I\alpha p}\,\pmb{v}^*_{\alpha I}\,,
\end{equation}
and
\begin{equation}
    p^t_{wp}(\pmb{x}_{w}) = \sum_{I=1}^{{N}_{n}} {S}_{Iwp}\,p^t_{wI}\,,
\end{equation}
where ${N}_{n}$ is the number of grid nodes that influence the material point, ${S}_{I\alpha p}$ is the shape function associated with node $I$ evaluated at the position of soil or material point $\alpha p$ ($\alpha = s\,\, or\,\, w$). The subscript $I$ represents the variables that are related to the grid nodes.

Using the same shape function, the test function $\delta\pmb{v}_{\alpha}$ and $\delta p$ can also be interpolated from their nodal values as
\begin{equation}\label{eq:test_func_v}
    \delta\pmb{v}_{\alpha} = \sum_{I=1}^{{N}_{n}} {S}_{I\alpha p}\,\delta\pmb{v}_{\alpha I} \,,
\end{equation}
and
\begin{equation}\label{eq:test_func_p}
    \delta{p} = \sum_{I=1}^{{N}_{n}} {S}_{I\alpha p}\,\delta{p}_{wI} \,.
\end{equation}

By treating the material points as quadrature points, we can solve the volume integral in the weak forms numerically by
\begin{equation}\label{eq:MPM_quadrature}
     \int_{\Omega_{\alpha}}^{} (\cdots) \, dV \approx  \sum_{\alpha p=1}^{{N}_{\alpha p}} (\cdots) \, V_{\alpha p}\,.
\end{equation}

The permanent information stored in the material points can be mapped to the grid. The mass and velocity of soil material points are mapped to the grid by
\begin{equation}
	\pmb{M}_{s I} = \sum_{s p=1}^{{N}_{s p}} {S}_{Is p} \, m_{s p} \,,
\end{equation}   
and
\begin{equation}
	\pmb{v}_{s I} = \frac{1}{\pmb{M}_{s I} } \sum_{s p=1}^{{N}_{s p}} {S}_{Is p} \, m_{s p} \, \pmb{v}_{s p} \,.
\end{equation}

Similarly, the porosity $n_{sp}$ and hydraulic conductivity $k_{sp}$ of soil are mapped to the grid by
\begin{equation}
	\pmb{n}_{I} = \frac{1}{\pmb{M}_{s I} } \sum_{sp=1}^{{N}_{sp}} {S}_{Isp} \, m_{s p} \, n_{sp} \,,
\end{equation}
and
\begin{equation}
	\pmb{k}_{I} = \frac{1}{\pmb{M}_{s I} } \sum_{sp=1}^{{N}_{sp}} {S}_{Isp} \, m_{s p} \, k_{sp} \,,
\end{equation}
where $n_{sp}$ and $k_{sp}$ are the porosity and hydraulic conductivity of soil, respectively, which are stored in the soil material points. Therefore, we can obtain them projected on the water material points by using
\begin{equation}
    n_{wp} = \sum_{I=1}^{{N}_{n}} {S}_{Iwp} \, \pmb{n}_{I} \,,
\end{equation}
and
\begin{equation}
    k_{wp} = \sum_{I=1}^{{N}_{n}} {S}_{Iwp} \, \pmb{k}_{I} \,,
\end{equation}
where $n_{wp}$ and $k_{wp}$ are the porosity and hydraulic conductivity of soil projected on the water material points, respectively. {To avoid the mapping error between soil and water material point, the grid water mass is obtained by using the grid porosity $\pmb{n}_{I}$ instead of $n_{wp}$:
\begin{equation}
	\pmb{M}_{w I} = \pmb{n}_{I}\sum_{w p=1}^{{N}_{w p}} {S}_{Iw p} \, \rho_w\, V_{w p} \,.
\end{equation}  
The velocity of water material points is mapped to the grid by
\begin{equation}\label{eq:v_wI}
	\pmb{v}_{w I} = \frac{\sum_{w p=1}^{{N}_{w p}} {S}_{I w p} \, (\rho_{w} \, V_{w p}^0) \, \pmb{v}_{w p}}{ \sum_{w p=1}^{{N}_{w p}} {S}_{I w p} \, (\rho_{w} \, V_{w p}^0) }  \,.
\end{equation} }

From the Eqs.~(\ref{eq:vI2vp})-(\ref{eq:v_wI}) and Eqs. (\ref{eq:msp})-(\ref{eq:mwp}), the weak forms in Section \ref{subsec:weak_form} can be discretised. The discretised momentum balance equations for the mixture and water are written as~\citep{Yuan2023}
\begin{equation}\label{eq:discretised_momentum_mixture}
	\pmb{M}_{sI}(\pmb{v}^*_{sI} - \pmb{v}^t_{sI}) + \pmb{M}_{wI}(\pmb{v}^*_{wI} - \pmb{v}^t_{wI}) = \Delta t(\pmb{F}^{ext}_I - \pmb{F}^{int}_I) \,,
\end{equation}
and
\begin{equation}\label{eq:discretised_momentum_water}
	\bar{\pmb{M}}_{wI} \pmb{v}^*_{wI} = \Delta t(\pmb{F}^{ext}_{wI} - \pmb{F}^{int}_{wI} + \pmb{Q}_I \pmb{v}^t_{sI}) + \pmb{M}_{wI} \pmb{v}^t_{wI} \,,
\end{equation}
where
%\begin{multline}
%	\pmb{F}^{ext}_{I} = \sum_{sp=1}^{{N}_{s p}} {S}_{Isp} \, m_{s p} \, \pmb{b} + \sum_{wp=1}^{{N}_{w p}} {S}_{Iwp} \, m_{w p} \, \pmb{b}\\
%	+ \sum_{sp=1}^{{N}_{s p}} \int_{\partial\Omega_{\tau}}^{} {S}_{Isp} \Bar{\pmb{\tau}} \, dS \,,
%\end{multline}
\begin{equation}
	\pmb{F}^{ext}_{I} = \sum_{sp=1}^{{N}_{s p}} {S}_{Isp} \, m_{s p} \, \pmb{b} + \sum_{wp=1}^{{N}_{w p}} {S}_{Iwp} \, m_{w p} \, \pmb{b}
	+ \sum_{sp=1}^{{N}_{s p}} \int_{\partial\Omega_{\tau}}^{} {S}_{Isp} \Bar{\pmb{\tau}} \, dS \,,
\end{equation}
\begin{equation}\label{eq:Fint}
	\pmb{F}^{int}_I = \sum_{sp=1}^{{N}_{s p}} \nabla {S}_{I sp} {\pmb{\sigma}'}^t_{sp}\,V_{sp}^t  + \sum_{wp=1}^{{N}_{w p}} \nabla {S}_{I wp}\, p^t_{wp}\,\pmb{I}\,V_{wp}^t \,,
\end{equation}
\begin{equation}
	\pmb{F}^{ext}_w = \sum_{wp=1}^{{N}_{w p}} {S}_{I wp} \, m_{w p} \, \pmb{b} \,,
\end{equation}
\begin{equation}\label{eq:Fint_w}
	\pmb{F}^{int}_w = \pmb{n}_I\sum_{wp=1}^{{N}_{w p}} \nabla {S}_{I wp}\, p^t_{wp}\,\pmb{I}\,V_{wp}^t \,,
\end{equation}
\begin{equation}
	\pmb{Q}_I = \frac{g\,\pmb{n}_I}{\pmb{k}_I}\pmb{M}_{wI}  \,,
\end{equation}
\begin{equation}
	\bar{\pmb{M}}_{wI} = \pmb{M}_{wI} + \pmb{Q}_I\Delta t \,.
\end{equation}
For every time step, the intermediate velocity of water $\pmb{v}^*_{wI}$ is obtained by solving Eq.~(\ref{eq:discretised_momentum_water}) explicitly. Then, the intermediate velocity of soil $\pmb{v}^*_{sI}$ can be calculated by solving Eq.~(\ref{eq:discretised_momentum_mixture}) also explicitly.

Following the same procedure, {by treating the material points as quadrature points (Eqs. (\ref{eq:test_func_v})-(\ref{eq:MPM_quadrature})) and moving the last two terms of the weak form (Eq. (\ref{eq:possion_eqn_weak})) to the right-hand-side}, the discretised Poisson equation is written as
\begin{equation}\label{eq:discretised_Poisson}
	\Delta t\,\pmb{L}\cdot(p^{t+1}_I - p^t_I) = \pmb{D}_{s}\cdot\pmb{v}^*_{sI} + \pmb{D}_{w}\cdot\pmb{v}^*_{wI} \,,
\end{equation}
where
%\begin{multline}
%	\pmb{L} =\frac{1-\pmb{n}_I}{\rho_s} \sum_{J=1}^{{N}_{n}} \sum_{sp=1}^{{N}_{sp}} V_{sp} \, \nabla {S}^p_{Isp} \, \nabla {S}^p_{Jsp} \\
%	+ \frac{\pmb{n}_I}{\rho_w}\sum_{J=1}^{{N}_{n}} \sum_{wp=1}^{{N}_{wp}} V_{wp} \, \nabla {S}^p_{Iwp} \, \nabla {S}^p_{Jwp} \,,
%\end{multline}
\begin{equation}
	\pmb{L} =\frac{1-\pmb{n}_I}{\rho_s} \sum_{J=1}^{{N}_{n}} \sum_{sp=1}^{{N}_{sp}} V_{sp}^t \, \nabla {S}^p_{Isp} \, \nabla {S}^p_{Jsp} 
	+ \frac{\pmb{n}_I}{\rho_w}\sum_{J=1}^{{N}_{n}} \sum_{wp=1}^{{N}_{wp}} V_{wp}^t \, \nabla {S}^p_{Iwp} \, \nabla {S}^p_{Jwp} \,,
\end{equation}
\begin{equation}
	\pmb{D}_{s} =(1-\pmb{n}_I)\sum_{J=1}^{{N}_{n}} \sum_{sp=1}^{{N}_{sp}} V_{sp}^t \, {S}^p_{Isp} \, \nabla {S}_{Jsp} \,,
\end{equation}
\begin{equation}
	\pmb{D}_{w} =\pmb{n}_I\sum_{J=1}^{{N}_{n}} \sum_{wp=1}^{{N}_{wp}} V_{wp}^t \, {S}^p_{Iwp}  \, \nabla {S}_{Jwp} \,.
\end{equation}
where ${S}^p_{I\alpha p}$ is the shape function for the water pressure field, which has only one degree of freedom because the water pressure is uniform in any direction. In contrast, the shape function for displacement field ${S}_{I\alpha p}$ has multi-degree of freedom. For simplicity, we use equal-order shape functions for pressure and displacement fields in this research. Given the intermediate velocities for soil and water phases $\pmb{v}^*_{\alpha I}$, the pressure increment $(p^{t+1}_I - p^t_I)$ can be obtained by solving Eq.~(\ref{eq:discretised_Poisson}) implicitly.

Finally, the second term of momentum balance equations for the mixture and water can be discretised as
\begin{equation}\label{eq:discretised_mixture_correction}
	\pmb{M}_{sI} (\pmb{v}^{t+1}_{sI} - \pmb{v}^{*}_{sI}) + \pmb{M}_{wI} (\pmb{v}^{t+1}_{wI} - \pmb{v}^{*}_{wI}) = \Delta t \, \pmb{G}(p^{t+1}_I - p^t_I) \,,
\end{equation}
and
\begin{equation}\label{eq:discretised_water_correction}
	 \pmb{M}_{wI} (\pmb{v}^{t+1}_{wI} - \pmb{v}^{*}_{wI}) = \Delta t \, \pmb{G}_w(p^{t+1}_I - p^t_I) \,,
\end{equation}
where
\begin{equation}
	\pmb{G} = \pmb{G}_s + \pmb{G}_w \,,
\end{equation}
\begin{equation}
	\pmb{G}_s = (1-\pmb{n}_I)\sum_{J=1}^{{N}_{n}} \sum_{sp=1}^{{N}_{sp}} V_{sp}^t \, {S}_{Isp} \nabla {S}^p_{Jsp} \,,
\end{equation}
\begin{equation}
	\pmb{G}_w = \pmb{n}_I\sum_{J=1}^{{N}_{n}} \sum_{wp=1}^{{N}_{wp}} V_{wp}^t \, {S}_{Iwp} \nabla {S}^p_{Jwp} \,.
\end{equation}

Using Eqs. (\ref{eq:discretised_mixture_correction}) and~(\ref{eq:discretised_water_correction}), the velocity of the soil and water phases can be updated.

\subsection{Constitutive models}
As we can see in Eq. (\ref{eq:Fint}), soil's effective stress ${\pmb{\sigma}'}^t_{sp}$ needs to be updated for each time step. 
In this research, a finite strain (or large strain) elastoplastic constitutive framework is used to update the stress. For a large deformation analysis, finite strain theory is a more rational choice \citep{Xie2023,Xie2022}. The interested reader may refer to the books written by \citet{Bathe1996} and \citet{deSouza2008} for more details of the finite strain constitutive model. 

In this coupled MPM formulation, the constitutive model is evaluated only at the soil material point level. To increase readability, we omit the subscript $sp$ for the notation in this section. Similarly, the stress in this section is always an effective stress, and its superscript is also omitted. 

In the finite strain theory, the left Cauchy-Green strain $\pmb{B}$ is defined as
\begin{equation}\label{eq:B}
	\pmb{B} = \pmb{F}\,\pmb{F}^T \,,
\end{equation}
and the left elastic Cauchy-Green strain ${\pmb{B}}^e$ is treated as a state variable, which keeps updated from time step $t$ to $t+1$ \citep{deSouza2008}:
\begin{equation}\label{eq:Be}
	\pmb{B}^{e,t+1} = \Delta\pmb{F}\, \pmb{B}^{e,t}\, \Delta\pmb{F}^T,
\end{equation}
where 
\begin{equation}
	\pmb{F} = \frac{\partial \pmb{x}}{\partial \pmb{X}} = \frac{\partial \left( \pmb{X}+\pmb{u} \right)}{\partial \pmb{X}} = \pmb{I} + \frac{\partial \pmb{u}}{\partial \pmb{X}}
\end{equation}
is the deformation gradient.

The logarithmic strain $\pmb{E}^e$ is one of the most commonly used finite strain measurements \citep{Xie2023b}, which is related to ${\pmb{B}}^e$ by
\begin{equation}\label{eq:log_strain}
	\pmb{E}^{e} = \frac{1}{2}\ln{\pmb{B}^{e}}.
\end{equation}
Note that the $\ln{(\Box)}$ operation in Eq.~(\ref{eq:log_strain}) requires eigendecomposition of ${\pmb{B}}^e$:
\begin{equation}
    \pmb{B}^e = \pmb{P}\,\pmb{D}^e\,\pmb{P}^T = \pmb{P}\,\text{diag}(\hat{\pmb{D}}^e)\,\pmb{P}^T \,,
\end{equation}
where $\pmb{P}$ is a rotation matrix whose columns are the normalised eigenvectors of $\pmb{B}^e$ and $\pmb{D}^e$ is a diagonal matrix whose diagonal entries $\hat{\pmb{D}}^e$ are the eigenvalues of $\pmb{B}^e$; the notation hat $\hat{(\Box)}$ here represents a column vector in principal space; the diagonal entries $\hat{\pmb{D}}^e$ are the principal stretches; the operation $\text{diag}(\Box)$ is to build a diagonal matrix whose only diagonal entries are non-zero. 
Therefore, Eq. (\ref{eq:log_strain}) yields
\begin{equation}\label{eq:log_strain2}
    \pmb{E}^{e} = \frac{1}{2}\,\pmb{P}\,\text{diag}(\ln{\hat{\pmb{D}}}^e)\,\pmb{P}^T \,.
\end{equation}
The operation $\ln{(\Box)}$ in Eq. (\ref{eq:log_strain2}) takes the logarithm of each component of the vector of $\hat{\pmb{D}}^e$.

Kirchhoff stress $\pmb{\tau}$ is usually used together with the logarithmic strain in a finite strain constitutive model. The Kirchhoff stress $\pmb{\tau}$ can be obtained by
\begin{equation}
    \pmb{\tau} = \mathbb{C}^{ep}:\pmb{E}^{e} \,,
\end{equation} 
where $\mathbb{C}^{ep}$ is the elastoplastic tangent operator which is a fourth-order tensor. The Cauchy stress $\pmb{\sigma}$ is related to Kirchhoff stress $\pmb{\tau}$ by
\begin{equation}
    \pmb{\sigma} = \frac{\pmb{\tau}}{J} \,,
\end{equation}
where $J = \det{\pmb{F}}$ is the Jacobian of deformation gradient, representing volume change. 

For an elastoplastic model, a yield criterion is required to decide whether the trail stress exceeds the limit, followed by a stress correction procedure if the stress actually has passed over the allowed value. In this research, the Mohr-Coulomb and the Nor-Sand models are used to represent the nonlinear behaviour of the soil. For the Mohr-Coulomb model, we strictly follow the algorithm documented by \citet{deSouza2008}. For the Nor-Sand model, originally proposed by \citet{Jefferies1993}, we follow its implicit return-mapping algorithm derived by \citet{Borja2006}. Using the Nor-Sand model for large deformation analysis requires some additional modifications. These modifications are included in Appendix~\ref{secA1} together with a validation.

\subsection{Stabilisation}
\subsubsection{Water phase}
Although the fractional step method satisfies the inf-sup condition~\citep{Kularathna2021}, additional stabilisation of the water pressure is still required for the incremental fractional step method. \citet{Yuan2023} introduced some artificial compressibilities to mitigate the oscillation of the water pressure for the incremental fractional step method. However, an additional artificial parameter is required in this approach. 

In contrast, a node-wise mapping and remapping technique is used in~\citet{Kularathna2021} to stabilise the water pressure. We adopt this approach because it is simple and efficient and no additional parameters are needed. This stabilisation method includes two steps: mapping the original water pressure at the material point level $p^t_{wp}$ to the related nodes and then mapping back to the governed material points. The following equation mathematically explains this procedure:
\begin{equation}
	\Bar{p}^t_{wp} = \sum_{I=1}^{{N}_{n}} {S}_{Iwp} \left(\frac{1}{m_{wI}}\sum_{p=1}^{{N}_{wp}} {S}_{Iwp} \, {m}_{wp} \, p^t_{wp}\right).
\end{equation}
where $m_{wI} = \sum_{wp=1}^{{N}_{wp}} {S}_{Iwp} \, m_{wp}$ is the nodal mass of water with a single degree of freedom, and $\Bar{p}^t_{wp}$ is the stabilised water pressure. In this research, we replace the original water pressure $p^t_{wp}$ in Eqs.~(\ref{eq:Fint}) and~(\ref{eq:Fint_w}) with the stabilised water pressure $\Bar{p}^t_{wp}$. 

\subsubsection{Soil phase}
%%%%%%%%%%
\begin{figure}
	\begin{center}
		\includegraphics[width=0.95\textwidth]{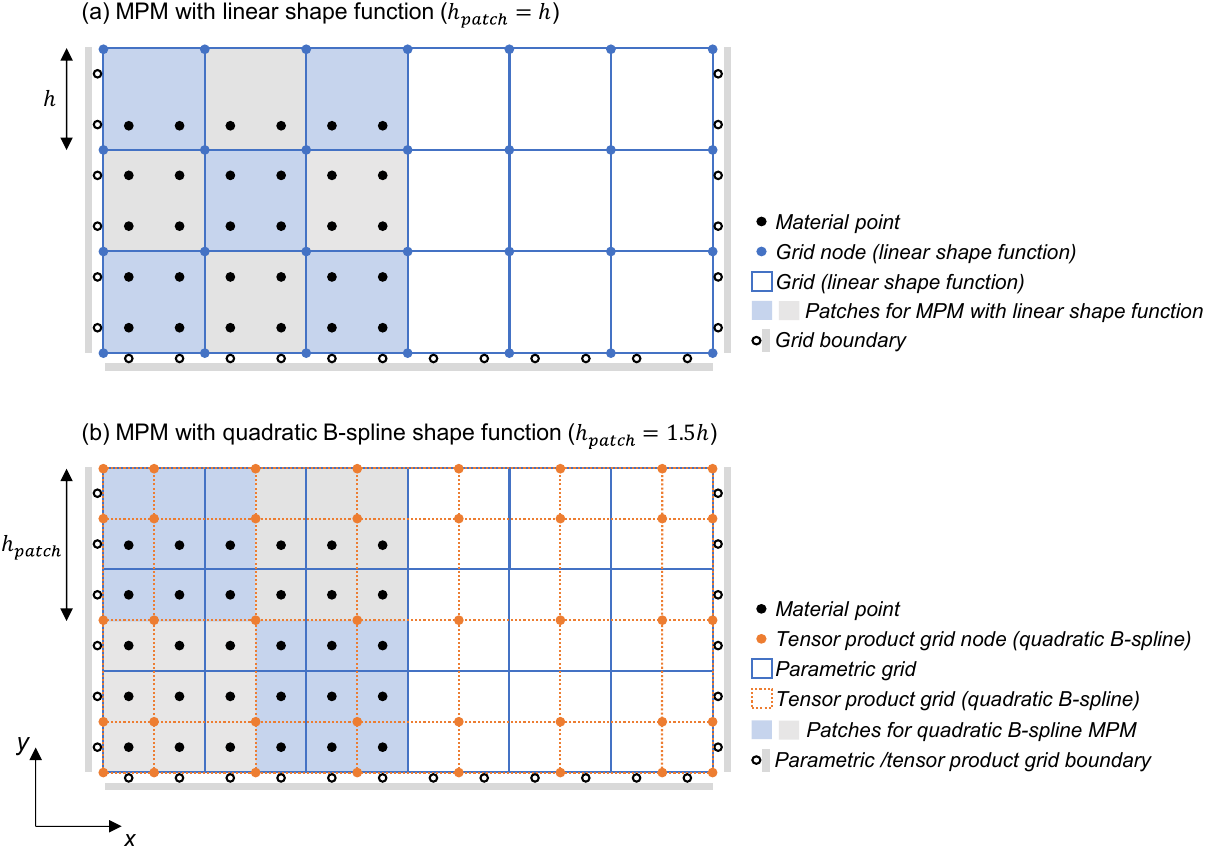}
		\caption{Illustration of the patch in the patch-wise averaging method. (\textbf{a}) Application on the MPM with linear shape function; (\textbf{b}) Application on the MPM with quadratic B-spline shape function}
		\label{mFbarPatch}
	\end{center}
\end{figure}
%%%%%%%%%%

As in the single-phase MPM, the volumetric locking instability also occurs in a coupled MPM. The F-bar method is one of the most commonly used stabilisation methods in MPM. The F-bar method was originally developed by \citet{deSouza1996} to counteract volumetric locking in FEM. \citet{Coombs2018} firstly applied the F-bar approach to stabilise MPM, which became a popular method to overcome this issue in MPM \citep{Moutsanidis2020, Wang2021, Telikicherla2022, Zhao2022, Sugai2023}. However, \citet{Xie2023b} discovered that the F-bar method cannot fully eliminate the stress oscillation in MPM under very large deformation. Therefore, \citet{Xie2023b} developed a modified F-bar approach to eliminate volumetric locking in MPM.

{In the F-bar method, $\Delta\pmb{F}$ in Eq.~(\ref{eq:Be}) is replaced with its stabilised form $\Delta\Bar{\pmb{F}}$ \citep{deSouza1996,Coombs2018}. Therefore, Eq.~(\ref{eq:Be}) yields to
\begin{equation}\label{eq:Fbar}
 	\pmb{B}^{e,t+1} = \Delta\Bar{\pmb{F}}{\pmb{B}}^{e,t} {\Delta\Bar{\pmb{F}}}^T \,.
\end{equation}
The stabilised form $\Delta\Bar{\pmb{F}}$ can be obtained by \citep{Coombs2018}
\begin{equation}
	\Delta\pmb{\Bar{F}} = \left(\frac{\Delta\Bar{J}}{\Delta{J}}\right)^{1/dof} \Delta\pmb{F},
\end{equation}
where $\Delta{J} = \det{(\Delta\pmb{F})}$, and \citep{Xie2023b}
\begin{equation}
	\Delta\Bar{J} = \sum_{I=1}^{{N}_{n}} {S}_{Isp} \left(\frac{1}{m_{sI}}\sum_{p=1}^{{N}_{sp}} {S}_{Isp} \, {m}_{sp} \, \Delta{J}\right).
\end{equation} } 

In the modified F-bar method, both $\Delta\pmb{F}$ and $\pmb{B}^{e,t}$ in Eq.~(\ref{eq:Be}) are replaced with their stabilised forms $\Delta\Bar{\pmb{F}}$ and $\Bar{\pmb{B}}^{e,t}$. As a result, Eq.~(\ref{eq:Be}) changes to
\begin{equation}\label{eq:mFbar}
 	\pmb{B}^{e,t+1} = \Delta\Bar{\pmb{F}}\Bar{\pmb{B}}^{e,t} {\Delta\Bar{\pmb{F}}}^T \,.
\end{equation}
In the original modified F-bar method \citep{Xie2023b}, the node-wise mapping and remapping technique is applied for both $\Delta\Bar{\pmb{F}}$ and $\Bar{\pmb{B}}^{e}$ to obtain their stabilised forms. {Although the node-wise stabilisation provides promising results in eliminating the volumetric locking, strain energy dissipates around the free surface of soil material points in this approach when it applies to stabilise $\Bar{\pmb{B}}^{e}$ \citep{Mast2012}. The original modified F-bar based on node-wise stabilisation of $\Bar{\pmb{B}}^{e}$ works on the Mohr-Coulomb model \citep{Xie2023b}. However, we found that applying node-wise mapping and remapping to $\pmb{B}^{e}$ will oversoften the structure on the Nor-Sand model. According to \citet{Mast2012}, a cell-wise averaging method, illustrated in Fig.~\ref{mFbarPatch}(a), is free of energy dissipation. This cell-wise method was designed for the MPM with linear shape function, such as general interpretation MPM (GIMPM). However, it is not feasible to apply a cell-wise averaging method to an MPM with a higher-order shape function, such as a B-spline MPM (BSMPM), because the shape function covers multiple cells in this case \citep{Xie2023b}. Therefore, in this research, a patch-wise averaging method is proposed to obtain the stabilised form $\Bar{\pmb{B}}^{e}$. The idea of the patch-wise averaging method is illustrated in Fig.~\ref{mFbarPatch}(b), retaining the same feature as the cell-wise method. The patch-wise averaging method covers multiple cells, resolving the challenge of the higher-order shape function. For a quadratic B-spline shape function used in this research, the patch size $h_{patch}$ is taken as 1.5 times as the cell size ($h$), as shown in Fig.~\ref{mFbarPatch}(b). For a cubic B-spline MPM, the same patch size (i.e. $h_{patch} = 1.5h$) can be used, resulting in a similar performance. Note that the patch is defined only once, and it remains unchanged during the simulation.} Applying the same concept of the F-bar method, the stabilised form $\Bar{\pmb{B}}^{e}$ in the modified F-bar method can be obtained by
\begin{equation}
	\Bar{\pmb{B}}^{e} = \bar{\pmb{F}}^e\,\bar{\pmb{F}}^{eT} =\left(\frac{\Bar{J}^e}{{J}^e}\right)^{2/3} \, {\pmb{F}}^e\,{\pmb{F}}^{eT} =  \left(\frac{\Bar{J}^e}{{J}^e}\right)^{2/3} \,\pmb{B}^{e},
\end{equation}
where ${J}^{e} = \det(\pmb{F}^{e}) = \sqrt{\det(\pmb{B}^{e})}$, and 
\begin{equation}
	\Bar{J}^{e}= \frac{\sum_{sp\in patch} {m}_{sp}\,{J}^{e}}{\sum_{sp\in patch} {m}_{sp}} \,.
\end{equation}

{The different performance of a cell-wise and patch-wise method on GIMPM and BSMPM is illustrated by a near incompressible elastic strip footing example in Appendix \ref{secA:mFbar}, indicating that the cell-wise method is only valid for a MPM with linear shape function. On the contrary, the patch-wise modified F-bar method effectively eliminates the volumetric locking in both GIMPM and BSMPM. Clearly, the resolution of the stress contour given by a patch-wise method depends on the mesh size as well. A contour with high resolution can be obtained by refining the mesh, as shown in Appendix \ref{secA:mFbar}. It is worth mentioning that the additional computational cost of these stabilisation methods is minimal. The accuracy of the analysis depends on the mesh size, but the localisation weakly depends on the resolution of the stress field. Therefore, it is not necessary to pursue an ultra-smooth stress contour by further refining the mesh if the problem is not mesh size sensitive.}

In addition to volumetric locking, MPM has stress oscillation around the physical body's boundary because of the ill-conditioned mass matrix. \citet{Coombs2023} solved this issue using the Ghost stabilisation technique. This method can only be applied to MPM with linear shape functions. \citet{Yamaguchi2021} developed an extended B-spline MPM to counteract this instability, which however, is usually ignored in explicit MPMs that adopt a lumped mass matrix because the inversion of the mass matrix is trivial in this case \citep{Coombs2023}. In this research, we also omit this instability due to computational cost reasons. Also, this instability occurs only around the boundary, and it is not severe compared with other issues such as the volumetric locking or the cell-crossing noise. 

\subsection{Semi-implicit two-phase double-point MPM algorithm}
In this research, we adopt the Modified Update-Stress-Last (MUSL) scheme for the algorithm and the FLuid-Implicit-Particle (FLIP) method to obtain the grid nodal velocities. See \citet{Zhang2016} for more details on MUSL and FLIP. The proposed semi-implicit two-phase double-point MPM algorithm, including stabilisation for the soil and water phases, is summarised as follows. 

\begin{enumerate}

\item[(1)] \emph{Shape functions and topology}

Obtain shape functions and their derivatives for the soil and water material points based on their coordinates: see \citet{Xie2023b} for the algorithm of the B-spline shape function.

\item[(2)] \emph{Map information from particles to nodes (P2N)}

(a). Obtain the nodal mass of the soil phase:
\begin{align*}
\pmb{M}_{s I} = \sum_{s p=1}^{{N}_{s p}} {S}_{I s p} \, m_{s p}
\end{align*}

(b). Obtain the nodal velocity of the soil phase:
\begin{align*}
\pmb{v}_{s I}^t = \frac{1}{\pmb{M}_{s I} } \sum_{s p=1}^{{N}_{s p}} {S}_{I s p} \, m_{s p} \, \pmb{v}_{s p}^t
\end{align*}

(c). Obtain the nodal porosity:
\begin{align*}
\pmb{n}_{I} = \frac{1}{\pmb{M}_{s I} } \sum_{sp=1}^{{N}_{sp}} {S}_{Isp} \, m_{s p} \, n_{sp}^t
\end{align*}

(d). Obtain the nodal hydraulic conductivity:
\begin{align*}
\pmb{k}_{I} = \frac{1}{\pmb{M}_{s I} } \sum_{sp=1}^{{N}_{sp}} {S}_{Isp} \, m_{s p} \, k_{sp}^t
\end{align*}

{(e). Obtain the nodal mass of the water phase:
\begin{align*}
	\pmb{M}_{w I} = \pmb{n}_{I}\sum_{w p=1}^{{N}_{w p}} {S}_{Iw p} \, \rho_w\, V_{w p}^t \,.
\end{align*} }

{(f). Obtain the nodal velocity of the water phase:
\begin{align*}
	\pmb{v}_{w I}^t = \frac{\sum_{w p=1}^{{N}_{w p}} {S}_{I w p} \, (\rho_{w} \, V_{w p}^0) \, \pmb{v}_{w p}^t }{ \sum_{w p=1}^{{N}_{w p}} {S}_{I w p} \, (\rho_{w} \, V_{w p}^0) } 
\end{align*} }

\item[(3)] \emph{Calculate the intermediate nodal velocity of the water phase explicitly}

(a). Obtain the internal nodal force of water phase:
\begin{align*}
\pmb{F}^{int}_w = \pmb{n}_I\sum_{wp=1}^{{N}_{w p}} \nabla {S}_{I wp}\, p^t_{wp}\,\pmb{I}\,V_{wp}^t
\end{align*}

(b). Obtain the external nodal force of water phase:
\begin{align*}
\pmb{F}^{ext}_w = \sum_{wp=1}^{{N}_{w p}} {S}_{I wp} \, m_{w p} \, \pmb{b}
\end{align*}

(c). Obtain the $\pmb{Q}_I$ and $\bar{\pmb{M}}_{wI}$ matrices:
\begin{align*}
	\pmb{Q}_I = \frac{g\,\pmb{n}_I}{\pmb{k}_I}\pmb{M}_{wI}
\end{align*}
\begin{align*}
	\bar{\pmb{M}}_{wI} = \pmb{M}_{wI} + \pmb{Q}_I\Delta t
\end{align*}

(d).  Obtain the intermediate nodal velocity of water phase:
\begin{align*}
     \pmb{v}^*_{wI} = \frac{\Delta t(\pmb{F}^{ext}_{wI} - \pmb{F}^{int}_{wI} + \pmb{Q}_I \pmb{v}^t_{sI}) + \pmb{M}_{wI} \pmb{v}^t_{wI}}{\bar{\pmb{M}}_{wI}}
\end{align*}

(e). Apply the Dirichlet boundary condition to $\pmb{v}^*_{wI}$.

\item[(4)] \emph{Calculate the intermediate nodal velocity of the soil phase explicitly}

(a). Obtain the internal nodal force of the mixture:
\begin{align*}
	\pmb{F}^{int}_I = \sum_{sp=1}^{{N}_{s p}} \nabla {S}_{I sp} {\pmb{\sigma}'}_{sp}^{t}\,V_{sp}^t  + \sum_{wp=1}^{{N}_{w p}} \nabla {S}_{I wp}\, p^t_{wp}\,\pmb{I}\,V_{wp}^t
\end{align*}

(b). Obtain the external nodal force of the mixture:
%\begin{multline*}
%	\pmb{F}^{ext}_{I} = \sum_{sp=1}^{{N}_{s p}} {S}_{Isp} \, m_{s p} \, \pmb{b} + \sum_{wp=1}^{{N}_{w p}} {S}_{Iwp} \, m_{w p} \, \pmb{b}\\
%	+ \sum_{sp=1}^{{N}_{s p}} \int_{\partial\Omega_{\tau}}^{} {S}_{Isp} \Bar{\pmb{\tau}} \, dS
%\end{multline*}
\begin{equation*}
	\pmb{F}^{ext}_{I} = \sum_{sp=1}^{{N}_{s p}} {S}_{Isp} \, m_{s p} \, \pmb{b} + \sum_{wp=1}^{{N}_{w p}} {S}_{Iwp} \, m_{w p} \, \pmb{b}
	+ \sum_{sp=1}^{{N}_{s p}} \int_{\partial\Omega_{\tau}}^{} {S}_{Isp} \Bar{\pmb{\tau}} \, dS
\end{equation*}

(c). Obtain the intermediate nodal velocity of soil phase:
\begin{equation*}
	\pmb{v}^*_{sI} = \frac{\Delta t(\pmb{F}^{ext}_I - \pmb{F}^{int}_I) - \pmb{M}_{wI}(\pmb{v}^*_{wI} - \pmb{v}^t_{wI})}{\pmb{M}_{sI}} + \pmb{v}^t_{sI}
\end{equation*}

(d). Apply the Dirichlet boundary condition to $\pmb{v}^*_{sI}$.

\item[(5)] \emph{Calculate the nodal pore pressure increment implicitly}

(a). Obtain the $\pmb{L}$, $\pmb{D}_{s}$, $\pmb{D}_{w}$ matrices:
%\begin{multline*}
%	\pmb{L} =\frac{1-\pmb{n}_I}{\rho_s} \sum_{J=1}^{{N}_{n}} \sum_{sp=1}^{{N}_{sp}} V_{sp} \, \nabla {S}^p_{Isp} \, \nabla {S}^p_{Jsp} \\
%	+ \frac{\pmb{n}_I}{\rho_w}\sum_{J=1}^{{N}_{n}} \sum_{wp=1}^{{N}_{wp}} V_{wp} \, \nabla {S}^p_{Iwp} \, \nabla {S}^p_{Jwp}
%\end{multline*}
\begin{equation*}
	\pmb{L} =\frac{1-\pmb{n}_I}{\rho_s} \sum_{J=1}^{{N}_{n}} \sum_{sp=1}^{{N}_{sp}} V_{sp}^t \, \nabla {S}^p_{Isp} \, \nabla {S}^p_{Jsp} 
	+ \frac{\pmb{n}_I}{\rho_w}\sum_{J=1}^{{N}_{n}} \sum_{wp=1}^{{N}_{wp}} V_{wp}^t \, \nabla {S}^p_{Iwp} \, \nabla {S}^p_{Jwp}
\end{equation*}

\begin{equation*}
	\pmb{D}_{s} =(1-\pmb{n}_I)\sum_{J=1}^{{N}_{n}} \sum_{sp=1}^{{N}_{sp}} V_{sp}^t \, {S}^p_{Isp} \, \nabla {S}_{Jsp}
\end{equation*}

\begin{equation*}
	\pmb{D}_{w} =\pmb{n}_I\sum_{J=1}^{{N}_{n}} \sum_{wp=1}^{{N}_{wp}} V_{wp}^t \, {S}^p_{Iwp}  \, \nabla {S}_{Jwp}
\end{equation*}

(b). Determine the free water surface: see \citet{Kularathna2017}.

(c). Apply the boundary condition for the free water surface, and solve nodal pore pressure increment: 
\begin{equation*}
	(p^{t+1}_I - p^t_I) = \pmb{L}^{-1}\cdot(\pmb{D}_{s}\cdot\pmb{v}^*_{sI} + \pmb{D}_{w}\cdot\pmb{v}^*_{wI})/\Delta t
\end{equation*}

\item[(6)] \emph{Update the nodal velocity of the soil and water phases}

(a). Obtain the $\pmb{G}_w$, $\pmb{G}_s$, and $\pmb{G}$ matrices:
\begin{equation*}
	\pmb{G}_w = \pmb{n}_I\sum_{J=1}^{{N}_{n}} \sum_{wp=1}^{{N}_{wp}} V_{wp}^t \, {S}_{Iwp} \nabla {S}^p_{Jwp}
\end{equation*}
\begin{equation*}
	\pmb{G}_s = (1-\pmb{n}_I)\sum_{J=1}^{{N}_{n}} \sum_{sp=1}^{{N}_{sp}} V_{sp}^t \, {S}_{Isp} \nabla {S}^p_{Jsp}
\end{equation*}
\begin{equation*}
	\pmb{G} = \pmb{G}_s + \pmb{G}_w
\end{equation*}

(b). Update the nodal velocity of the water phase:
\begin{equation*}
	 \pmb{v}^{t+1}_{wI} = \frac{\Delta t \, \pmb{G}_w(p^{t+1}_I - p^t_I)}{\pmb{M}_{wI}} + \pmb{v}^{*}_{wI} 
\end{equation*}

(c). Update the nodal velocity of the soil phase:
\begin{equation*}
	\pmb{v}^{t+1}_{sI} = \frac{\Delta t \, \pmb{G}(p^{t+1}_I - p^t_I) - \pmb{M}_{wI} (\pmb{v}^{t+1}_{wI} - \pmb{v}^{*}_{wI})}{\pmb{M}_{sI}}  + \pmb{v}^{*}_{sI}  
\end{equation*}

(d). Apply the Dirichlet boundary condition to $\pmb{v}^*_{wI}$ and $\pmb{v}^*_{sI}$.

\item[(7)] \emph{Map information from nodes to particles(N2P)}

(a). Update the coordinates of soil and water material points:
\begin{equation*}
    \pmb{x}_{\alpha p}^{t+1} = \pmb{x}_{\alpha p}^{t} + \Delta t \sum_{I=1}^{{N}_{n}} {S}_{I\alpha p} \, \pmb{v}_{\alpha I}^{t+1}
\end{equation*}

(b). Update the velocity of soil and water material points:
\begin{equation*}
    \pmb{v}_{\alpha p}^{t+1} = \pmb{v}_{\alpha p}^{t} + \sum_{I=1}^{{N}_{n}} {S}_{I\alpha p} \, (\pmb{v}_{\alpha I}^{t+1}-\pmb{v}_{\alpha I}^{t})
\end{equation*}

(c). Update the pore pressure of water material points:
\begin{equation*}
    {p}_{wp}^{t+1} = \sum_{I=1}^{{N}_{n}} {S}_{Iwp} \, p^{t+1}_I
\end{equation*}

\item[(8)] \emph{Update information at the material point level}

(a). Obtain the nodal displacement increment of soil phase:
\begin{equation*}
	\Delta\pmb{u}_{sI}^{t+1} = \frac{\Delta t}{\pmb{M}_{s I} } \sum_{s p=1}^{{N}_{sp}} {S}_{Isp} \, m_{sp} \, \pmb{v}_{sp}^{t+1}
\end{equation*}

(b). Obtain the deformation gradient increment of soil material points:
\begin{equation*}
    \Delta\pmb{F}_{sp}^{t+1} = \pmb{I} + \sum_{I=1}^{{N}_{n}} \, \nabla S_{Isp} \otimes \Delta\pmb{u}_{sI}^{t+1}
\end{equation*}

(c). Obtain the Jacobian of deformation gradient increment of soil material points:
\begin{equation*}
    \Delta J_{sp}^{t+1} = \det{\Delta\pmb{F}_{sp}^{t+1}}
\end{equation*}

(d). Update the volume of soil material points:
\begin{equation*}
    V_{sp}^{t+1} = \Delta J_{sp}^{t+1} \, V_{sp}^{t}
\end{equation*}

(e). Obtain the Jacobian of deformation gradient of soil material points:
\begin{equation*}
    J_{sp}^{t+1} = V_{sp}^{t+1}/V_{sp}^{0}
\end{equation*}

(f). Update the porosity of soil material points:
\begin{equation*}
    n_{sp}^{t+1} = 1-(1-n_{sp}^{0})/J_{sp}^{t+1}
\end{equation*}

(g). Update the hydraulic conductivity of soil material points based on Kozeny-Carman equation(see \citet{Carrier2003,Bandara2015} for more details):
\begin{equation*}
    k_{sp}^{t+1} = C_1\frac{(n_{sp}^{t+1})^3}{(1-n_{sp}^{t+1})^2}
\end{equation*}

(h). Obtain the porosity of soil projected on water material points:
\begin{equation*}
    n_{wp}^{t+1} = \sum_{I=1}^{{N}_{n}} {S}_{Iwp} \left(\frac{1}{\pmb{M}_{wI}}\sum_{sp=1}^{{N}_{sp}} {S}_{Isp} \, m_{sp} \, n_{sp}^{t+1} \right)
\end{equation*}

(I). Update the volume of water material points:
\begin{equation*}
    V_{wp}^{t+1} = n_{wp}^{0}/n_{wp}^{t+1} \, V_{wp}^{0} 
\end{equation*}

\item[(9)] \emph{Stabilise the water and soil phases}

(a). Stabilise the pore pressure:
\begin{equation*}
	{p}_{wp}^{t+1} = \sum_{I=1}^{{N}_{n}} {S}_{Iwp} \left(\frac{1}{m_{wI}}\sum_{p=1}^{{N}_{wp}} {S}_{Iwp} \, {m}_{wp} \, p_{wp}^{t+1}\right)
\end{equation*}

(b). Obtain the stabilised Jacobian of deformation gradient of soil material points:
\begin{equation*}
	\Bar{J}_{sp}^{t+1} = \sum_{I=1}^{{N}_{n}} {S}_{Isp} \left(\frac{1}{m_{sI}}\sum_{p=1}^{{N}_{sp}} {S}_{Isp} \, {m}_{sp} \, {J}_{sp}^{t+1}\right)
\end{equation*}

(c). Obtain the stabilised Jacobian of deformation gradient increment of soil material points:
\begin{equation*}
	\Delta\Bar{J}_{sp}^{t+1} = \sum_{I=1}^{{N}_{n}} {S}_{Isp} \left(\frac{1}{m_{sI}}\sum_{p=1}^{{N}_{sp}} {S}_{Isp} \, {m}_{sp} \, \Delta{J}_{sp}^{t+1}\right)
\end{equation*}

(d). Obtain the stabilised deformation gradient increment of soil material points:
\begin{equation*}
	\Delta\pmb{\Bar{F}}_{sp}^{t+1} = \left(\frac{\Delta\Bar{J}_{sp}^{t+1}}{\Delta{J}_{sp}^{t+1}}\right)^{1/dof} \Delta\pmb{F}_{sp}^{t+1},
\end{equation*}

(e). Obtain the Jacobian of elastic deformation gradient of soil material points:
\begin{equation*}
	{J}_{sp}^{e,t}=\sqrt{\det(\pmb{B}^{e,t}_{sp})}
\end{equation*}

{(f). Obtain the stabilised Jacobian of elastic deformation gradient of soil material points:
\begin{equation*}
	\Bar{J}_{sp}^{e,t}= \frac{\sum_{sp\in patch} {m}_{sp}\,{J}_{sp}^{e,t}}{\sum_{sp\in patch} {m}_{sp}}
\end{equation*} }

(g). Obtain the stabilised last converged elastic left Cauchy-
Green strain:
\begin{equation*}
	\Bar{\pmb{B}}^{e,t}_{sp} =  \left(\frac{\Bar{J}_{sp}^{e,t}}{{J}_{sp}^{e,t}}\right)^{2/3} \,\pmb{B}^{e,t}_{sp},
\end{equation*}

\item[(10)] \emph{Update stress and strain}

(a). Obtain the elastic trial left Cauchy-Green strain:
 \begin{equation*}
 	\pmb{B}^{e,trial}_{sp} = \Delta\Bar{\pmb{F}}^{t+1}_{sp}\Bar{\pmb{B}}^{e,t}_{sp} {\Delta\Bar{\pmb{F}}^{t+1}_{sp}}^T
 \end{equation*}

(b). Perform eigendecomposition to elastic trial left Cauchy-Green strain:
\begin{equation*}
    \pmb{B}^{e,trial}_{sp} = \pmb{P}\,\pmb{D}^e\,\pmb{P}^T = \pmb{P}\,\text{diag}(\hat{\pmb{D}}^e)\,\pmb{P}^T
\end{equation*}

(c). Update the Kirchhoff stress $\pmb{\tau}^{t+1}_{sp}$ and the elastic left Cauchy-Green strain $\pmb{B}^{e,t+1}_{sp}$ based on the constitutive model: see \citet{deSouza2008} for the algorithms of the Mohr-Coulomb model; see \citet{Borja2006} for the algorithms of the Nor-Sand model and the Appendix \ref{secA1} for its modification for a large deformation problem.

(d). Update the Cauchy stress:
\begin{equation*}
    {\pmb{\sigma}'}^{t+1}_{sp} = \frac{\pmb{\tau}^{t+1}_{sp}}{{J}_{sp}^{t+1}}
\end{equation*}

\end{enumerate}

\section{Numerical examples}\label{sec5}
In this section, some numerical examples are presented. First, the proposed two-phase double-point algorithm is validated with Terzaghi’s one-dimensional consolidation theory \citep{Terzaghi1943}. Then, the performance of this algorithm under large deformation is validated with {the analytical solution for one-dimensional large deformation consolidation derived by \citet{Xie2004} and }a Mohr-Coulomb slope stability problem performed by \citet{Yuan2023}. Finally, a landslide numerical simulation with Nor-Sand constitutive model is performed using both single-point and double-point formulations, to compare the results obtained with both approaches. The most widely used general interpretation MPM (GIMPM) cannot fully remove the cell-crossing instability, and we find that the convergence of the Nor-Sand constitutive model becomes challenging due to this instability. Therefore, these numerical examples are conducted using a quadratic B-spline shape function (see \citet{Xie2023b} for more details on B-spline MPM), to fully eliminate cell-crossing noise at an affordable computational cost. Although the cubic B-spline shape function is more accurate, a significantly lower computational cost is required for a quadratic shape function. Additionally, the convergence issue of the Nor-Sand model worsens if either the water or the soil phases are not stabilised. 

\subsection{One-dimensional consolidation}
\begin{figure}
	\begin{center}
		\includegraphics[width=0.45\textwidth]{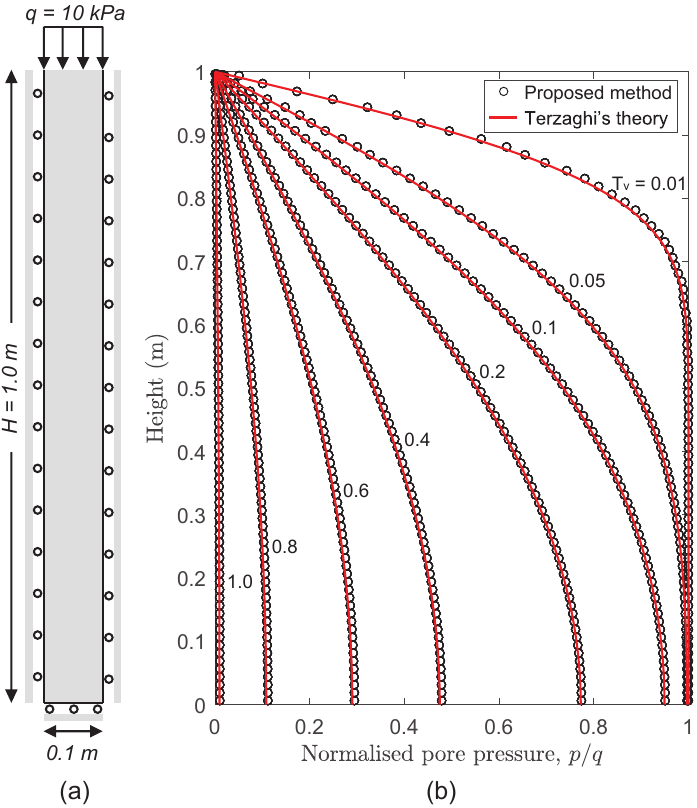}
		\caption{One-dimensional consolidation test. (\textbf{a}) Graphical illustration of the test; (\textbf{b}) Comparison of the pore pressure of the proposed method with Terzaghi’s solution}
		\label{1Dconsolidation}
	\end{center}
\end{figure}

Terzaghi’s one-dimensional consolidation is one of the traditional benchmark problems to validate the soil-water coupled numerical models \citep{Bandara2015, Navas2016, Yamaguchi2020, Kularathna2021, Pan2021, Zheng2022, Yuan2023}. The graphical illustration of this problem is shown in Fig. \ref{1Dconsolidation}(a). A 2-D plain strain soil column is modelled with a height of $1\,m$ and a width of $0.1\,m$. The roller boundary condition is applied on the bottom and lateral surfaces of this column (that is, at the boundary nodes). The impermeable boundary condition is automatically satisfied at these boundaries, where the displacement is constrained. No displacement is prevented at the top surface, where the pore water pressure is free to dissipate. A uniformly distributed load (UDL), $q = 10\,kPa$, is initially applied to its upper surface and held until the last step, and the gravitational force is ignored for this simulation. The soil column is assumed to be linear elastic, and the material parameters are summarised in Table. \ref{tab1D}. For this problem, a uniformly distributed $0.025\,m$ wide square grid with $2^2$ material points per cell is used, resulting in 640 material points in total for both water and soil phases. The time step for this problem is $10^{-4}\,s$, and the total time is about $1\,s$.
\begin{table}[h]
	\caption{Parameters for Terzaghi’s one-dimensional consolidation}\label{tab1D}%
	\begin{tabular}{@{}ll@{}}
		\toprule
		Parameter & Value  \\
		\midrule
		Water density $\rho_w [kg/m^3]$    & 1000 \\
		Soil grain density $\rho_s [kg/m^3]$    & 2650 \\
		Young’s modulus $E[MPa]$    & 10 \\
		Poisson’s ratio $\nu$    & 0 \\
		Initial porosity $n_0$    & 0.3 \\
		Initial hydraulic conductivity $k_0 [m/s]$    & $10^{-3}$ \\
		Gravitational acceleration $g[m/s^2]$    & 0 \\
		\bottomrule
	\end{tabular}
\end{table}
%%%%%%%%%%

According to Terzaghi’s one-dimensional consolidation theory \citep{Terzaghi1943}, the analytical pore pressure along the soil column can be obtained by
\begin{equation}
    p(t,z) = \frac{4q}{\pi} \sum_{i=1}^{\infty} \frac{1}{2i+1} \sin\left[\frac{(2i+1)\pi z}{2H}\right] e^{-(i+0.5)^2{\pi}^2T_v}\,,
\end{equation}
where
\begin{equation}
    T_v = \frac{k\,t}{\rho_w\,g\,m_v\,H^2}
\end{equation}
is the normalised time factor, 
\begin{equation}
    m_v = \frac{(1+\nu)(1-2\nu)}{(1-\nu)E},
\end{equation}
is the inverse of the pressure wave modulus or the constrained modulus, and $H$ is the height of the soil column.

Normalised pore pressures along the soil column given by both numerical simulation and analytical solution are plotted in Fig. \ref{1Dconsolidation}(b) with different normalised time factors, showing that the proposed two-phase, double-point approach is in excellent agreement with Terzaghi’s analytical solutions.

\subsection{One-dimensional large deformation consolidation}
{The proposed method has been validated in the small deformation regime in the previous numerical example. In this example, the proposed method is tested under the large deformation by comparing it with the analytical solution for the one-dimensional large deformation consolidation derived by \citet{Xie2004}. This analytical solution has also been used to validate many coupled MPMs \citep{Zheng2021, Martinelli2022,Yuan2023}. The setup for the numerical model is the same as Terzaghi’s one-dimensional consolidation, except Young's modulus of the soil column $E$ has been reduced to 1 MPa and the applied UDL $q$ has been increased to 200 kPa to ensure a large deformation during the consolidation. }

{According to the one-dimensional large deformation consolidation theory \citet{Xie2004} proposed, the settlement on the top surface of the soil column can be calculated by
\begin{equation}
     S_{t} = S(t,\,z = 0) = H\left(1 - e^{-m_v\,q}\right) \left(1 - \sum_{m=1}^{\infty} \frac{2}{{M}^2} \, e^{-4M^{2}\,T_v}\right) \,,
\end{equation}
where $M = (m-\frac{1}{2})\pi \,, \,\,\,\, m = 1,\, 2,\, 3,\,\cdots$.
Therefore, the ultimate settlement is:
\begin{equation}
    S_{\infty} = S(t = \infty,\,z = 0) = H(1 - e^{-m_v\,q}) \,.
\end{equation}
The degree of consolidation $U_s$ can be defined as the ratio between the top surface of the soil column $S_t$ and the ultimate settlement $S_{\infty}$: 
\begin{equation}
    U_s = \frac{S_t}{S_{\infty}} = 1 - \sum_{m=1}^{\infty} \frac{2}{{M}^2} \, e^{-4M^{2}\,T_v} \,.
\end{equation}
The pore pressure along the soil column can be obtained by \citep{Xie2004}
\begin{equation}
    p(t,\,z) = \frac{1}{m_v} \ln{\left[1+(e^{m_v \, q}-1)\sum_{m=1}^{\infty} \frac{2}{M} \, \sin{\left(\frac{2Mz}{H}\right)} \, e^{-4M^{2}\,T_v}\right]}  \,.
\end{equation} }

{Same as \citet{Zheng2021, Martinelli2022,Yuan2023}, $m_v$ and $T_v$ are obtained using the one defined in the previous Terzaghi’s one-dimensional consolidation test. The soil column is assumed to be linear elastic with the material parameters summarised in Table. \ref{tab1Dlarge}. For this problem, a uniformly distributed 0.05 wide square grid with $2^2$ material points per cell is used, resulting in 160 material points in total for both water and soil phases. The time step for this problem is $5\times10^{-4}$ s, and the total time is about 20 s to reach $T_v = 2$ s. For a comparison purpose, the simulations have been carried out using the proposed method with the constitutive models based on small and large strain frameworks, respectively. The interested reader can refer to \citet{deSouza2008, Xie2023} for the difference between the small and large strain constitutive models. }
\begin{table}[h]
	\caption{Parameters for one-dimensional large deformation consolidation}\label{tab1Dlarge}%
	\begin{tabular}{@{}ll@{}}
		\toprule
		Parameter & Value  \\
		\midrule
		Water density $\rho_w [kg/m^3]$    & 1000 \\
		Soil grain density $\rho_s [kg/m^3]$    & 2650 \\
		Young’s modulus $E[MPa]$    & 1 \\
		Poisson’s ratio $\nu$    & 0 \\
		Initial porosity $n_0$    & 0.3 \\
		Initial hydraulic conductivity $k_0 [m/s]$    & $10^{-3}$ \\
		Gravitational acceleration $g[m/s^2]$    & 0 \\
		\bottomrule
	\end{tabular}
\end{table}

\begin{figure}
	\begin{center}
		\includegraphics[width=0.95\textwidth]{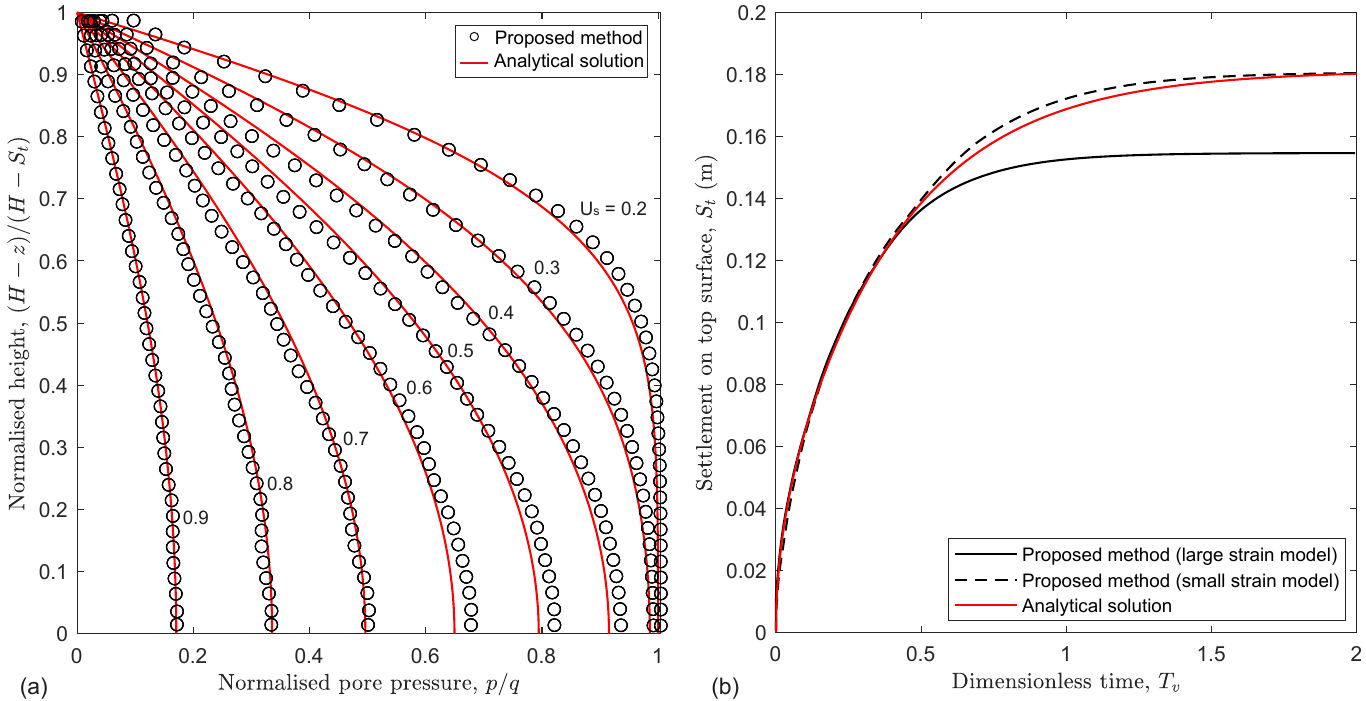}
		\caption{One-dimensional large deformation consolidation test. (\textbf{a}) Comparison of the pore pressure of the proposed method with the analytical solution at different degrees of consolidation $U_s$; (\textbf{b}) Comparison of the settlement $S_t$ of the proposed method with analytical solution}
		\label{1Dlarge}
	\end{center}
\end{figure}

{The pore pressure distribution at different degrees of consolidation and the time history of settlement given by the proposed method are compared with the analytical solution in Fig. \ref{1Dlarge}. As shown in Fig. \ref{1Dlarge} (a), the pore pressure $p$ is normalised by the applied UDL $q$, and the vertical position of the soil column $H-z$ is normalised by the current height of the column $H-S_t$. The proposed method has an excellent match with the analytical solution in terms of pore water pressure dissipation under the large deformation. The normalised pore pressure distributions given by the small and large strain constitutive models are the same, which is represented by the same group of data points named as the proposed method in Fig. \ref{1Dlarge} (a). However, their settlement time histories are quite different, as shown in Fig. \ref{1Dlarge} (b). At the relatively small strain range, both small and large strain constitutive models march the analytical solution perfectly. However, as the soil column undergoes large deformation, the small strain constitutive model overestimates the settlement by about 17\%. A similar amount of overestimation can also be observed in the previous studies \citep{Li2004,Navas2017,Navas2018} that compare the small and large strain constitutive model using a similar large deformation consolidation problem. As we can see in  Fig. \ref{1Dlarge} (b), the result given by the small strain constitutive model matches the analytical solution proposed by \citet{Xie2004} very well, implying an updated analytical solution may be needed for a large strain consolidation problem. Nevertheless, the performance of the proposed numerical model under large deformation can be validated. }
%%%%%%%%%%

\subsection{Slope stability with Mohr-Coulomb model}
\begin{figure}
	\begin{center}
		\includegraphics[width=0.45\textwidth]{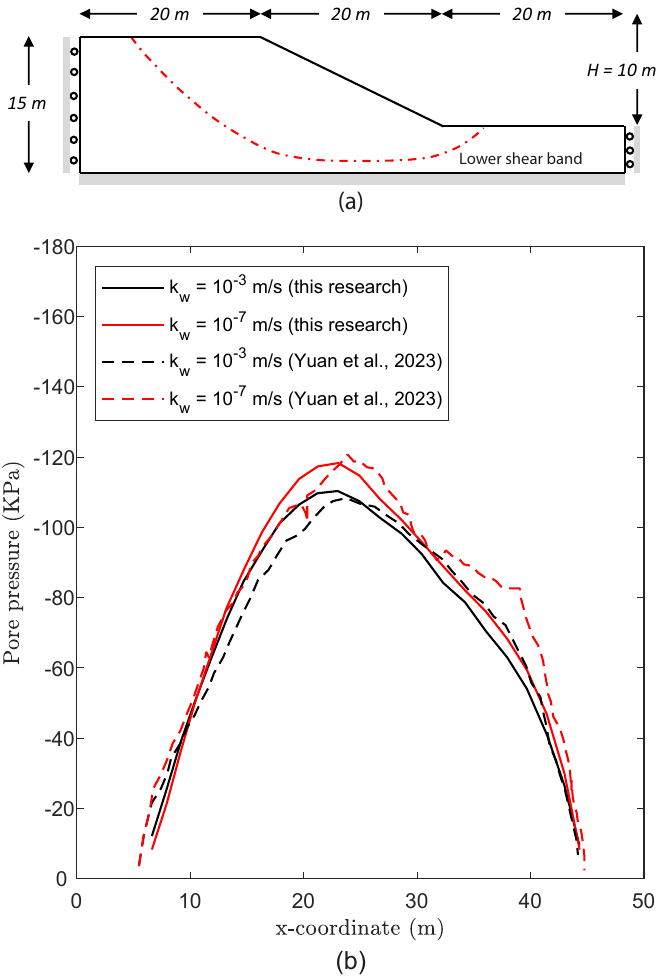}
		\caption{Slope stability with Mohr-Coulomb model (\textbf{a}) Graphical illustration of the numerical model; (\textbf{b}) Comparison of pore pressure along the lower shear band at t = 3.0 s from the proposed method with \citet{Yuan2023}}
		\label{MCslope}
	\end{center}
\end{figure}

In this section, the Mohr-Coulomb slope stability problem studied by \citet{Yuan2023} is used for validation of the methodology proposed in this paper. Fig. \ref{MCslope}(a) shows the geometry and boundary conditions for the numerical model. The roller and fixed boundary conditions are applied on the sides and bottom of the model, respectively. {Impermeable boundaries are automatically formed where those roller and fixed boundary conditions are applied. No pressure boundary condition is applied on the water phase, and the pore pressure is allowed to dissipate freely on the free water surface.} The material parameters for the Mohr-Coulomb model are summarised in Table \ref{tabMC}. A strength reduction factor (SRF) is applied to both soil cohesion $c$ and friction angle $\phi$ to trigger the failure. Therefore, the reduced cohesion of the soil $c_f$ and friction angle $\phi_f$ can be obtained by \citep{Yuan2023}
\begin{equation}
    c_f = c/SRF \,\,\,\, and \,\,\,\, \phi_f = \arctan{(\tan{\phi}/SRF)}\,.
\end{equation}
In this study, $SRF = 2.0$ is applied to compare the results given by \citet{Yuan2023}.
\begin{table}[h]
	\caption{Parameters for slope stability with Mohr-Coulomb model}\label{tabMC}%
	\begin{tabular}{@{}ll@{}}
		\toprule
		Parameter &Value  \\
		\midrule
		Water density $\rho_w [kg/m^3]$    & 1000 \\
		Soil grain density $\rho_s [kg/m^3]$    & 2500 \\
		Young’s modulus $E[MPa]$    & 100 \\
		Poisson’s ratio $\nu$    & 0.3 \\
		Cohesion $c[kPa]$    & 20 \\
		Friction angle $\phi[^{\circ}]$    & 30 \\
		Dilation angle $\psi[^{\circ}]$    & 0 \\
		Strength reduction factor $SRF$    & 2.0 \\
		Initial porosity $n_0$    & 0.4 \\
		Initial hydraulic conductivity $k_0 [m/s]$    & $10^{-3}$ \& $10^{-7}$ \\
		Gravitational acceleration $g[m/s^2]$    & -9.81 \\
		\bottomrule
	\end{tabular}
\end{table}
%%%%%%%%%%

For this numerical example, a uniformly distributed $0.5\,m$ wide square grid with $4^2$ material points per cell is used, resulting in 38,320 material points in total for both water and soil phases. The mesh size remains the same as \citet{Yuan2023}, and $4^2$ material points per cell are sufficient to ensure the accuracy of this problem. The time step is $0.002\,s$, which is about 0.92 CFL condition. The final time is set as 6 s, and at this time the slope is almost static. The initialisation and loading procedures are the same as the ones proposed by \citet{Yuan2023}. {It is worth mentioning that the same time step is used for both cases with high and low hydraulic conductivities. In other words, the critical time step of the proposed formulation is only related to the CFL condition for the single-phase MPM. This is not the case in the semi-implicit MPM proposed by \cite{Yamaguchi2020}, leading to a very small time step in the low hydraulic conductivity case.} 

Fig. \ref{MCslope}(b) compares the results of the pore pressure along the lower shear band at 3 s obtained by the proposed method and \citet{Yuan2023}. The position of the lower shear band, defined by \citet{Yuan2023}, is sketched in Fig. \ref{MCslope}(a). In general, the pore pressures given by the proposed method agree with \citet{Yuan2023} for both low and high hydraulic conductivity cases, indicating the validation of the proposed method in this case in which large deformation occurs. The discrepancies between the present research and \citet{Yuan2023} may be due to the different shape functions. The GIMPM is used in \citet{Yuan2023}, and the higher order B-spline MPM is used in the current research. Therefore, the positions of the shear band may be slightly different in both cases. However, higher-order MPM usually yields more accurate results \citep{Xie2023}. The different stabilisation methods may also be one of the reasons. Nevertheless, a smoother pore pressure distribution can be observed in the curves given by this research, as shown in Fig. \ref{MCslope}(b).
\begin{figure}
	\begin{center}
		\includegraphics[width=0.9\textwidth]{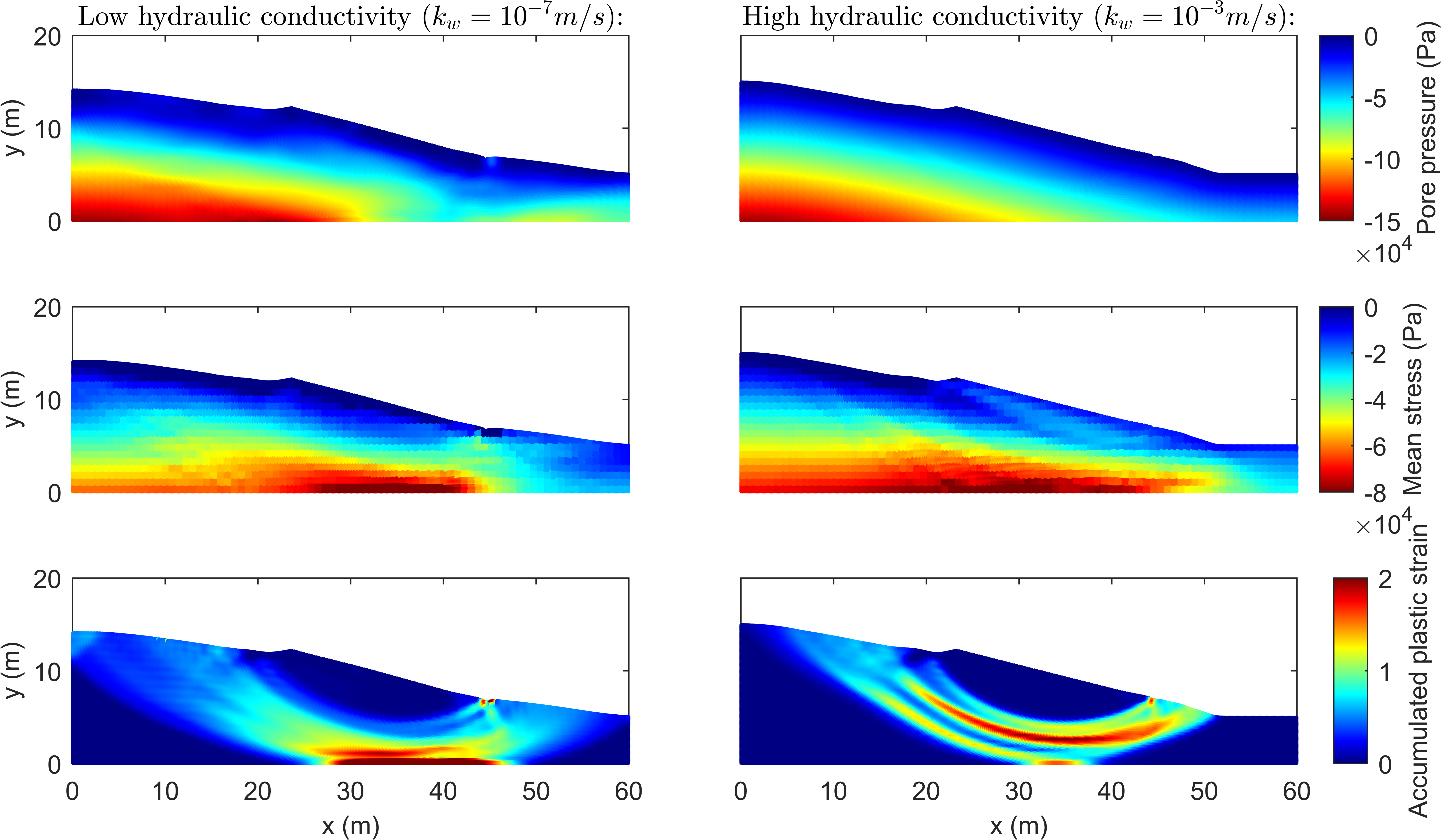}
		\caption{Contours of the slope stability with Mohr-Coulomb model at the final step (t = 6s). (\textbf{a}) the case of low hydraulic conductivity ($k_w = 10^{-7}$); (\textbf{b}) the case of high hydraulic conductivity ($k_w = 10^{-3}$)}
		\label{MCslope_result}
	\end{center}
\end{figure}

The pore pressures, the accumulation of plastic strain, and the mean effective stress contours for the case of low and high hydraulic conductivity in the final step are represented in the left and right parts of Fig. \ref{MCslope_result}, respectively. The pore pressure contours are plotted using the water material points. The accumulated plastic strain {obtained from the return mapping algorithm \citep{deSouza2008},} and the mean effective stress contours are represented by using the soil material points. However, for this numerical problem, the relative displacement between the water and soil material point is very small even in the case of high hydraulic conductivity, being 0.0014 m and 0.0065 m, respectively. A large relative displacement between water and soil material points can be expected in a dilative failure, which will be shown in the next numerical example. 

\begin{figure}
	\begin{center}
		\includegraphics[width=0.95\textwidth]{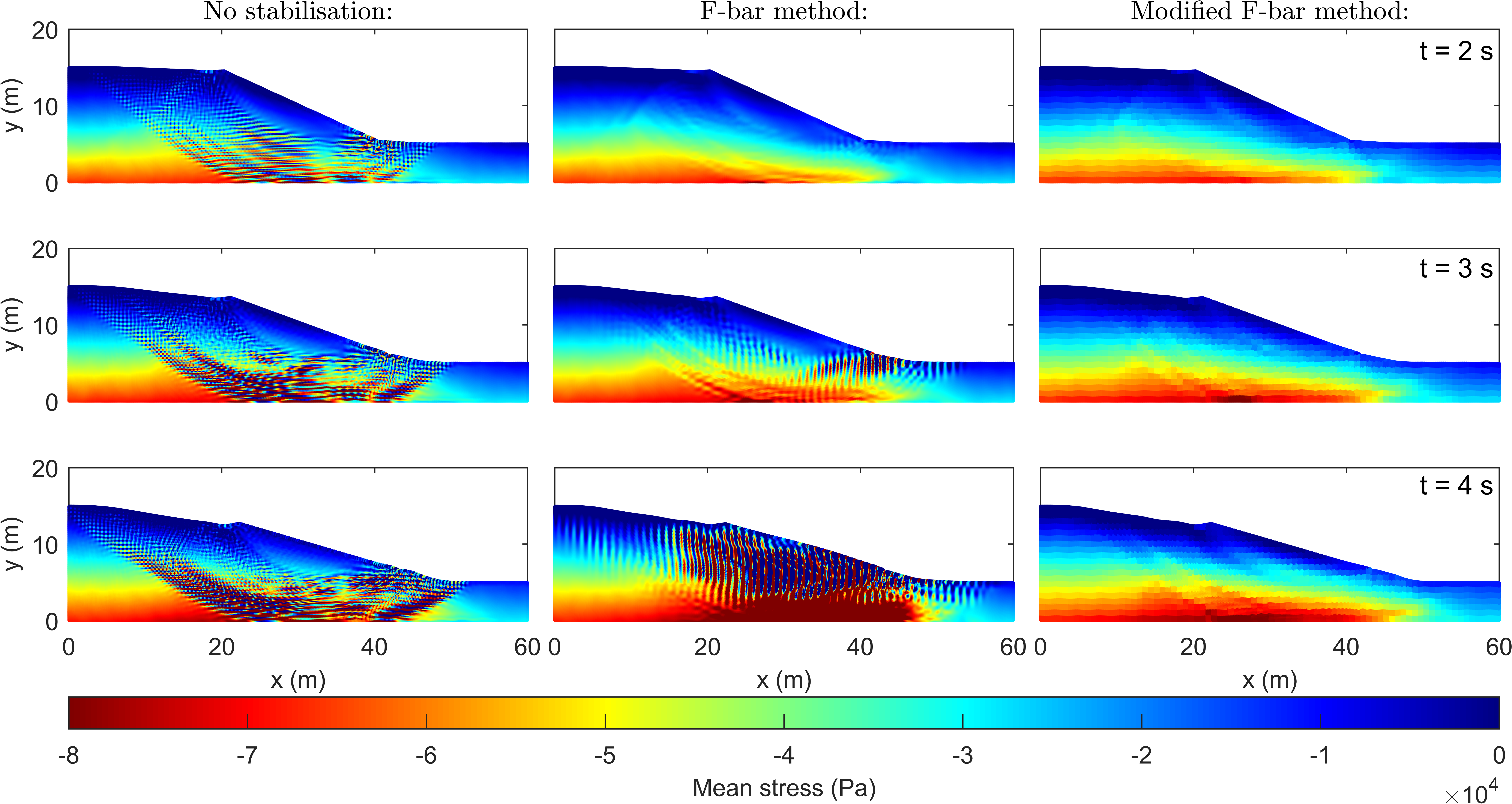}
		\caption{Soil mean effective stress contours of the slope stability with Mohr-Coulomb model under high hydraulic conductivities with different stabilisation methods at various time snapshots.}
		\label{2MP_MCslope_compareStab}
	\end{center}
\end{figure}

As we can see in Fig. \ref{MCslope_result}, the pore pressure contours are very smooth in both cases, and they have very similar distributions compared with \citet{Yuan2023}. The mean effective stress contours have relatively low resolution compared with the pore pressure contours. This is because the cell-wise averaging procedure is adopted in the modified F-bar method, while the node-wise averaging procedure is adopted for the pore pressure. However, the stress oscillation due to volumetric locking is significantly removed, and the resolution of the stress contours can be improved by using a finer mesh {as demonstrated in Appendix \ref{secA:mFbar}}. Unfortunately, it is not possible to compare our stress contours with \citet{Yuan2023} because the one under large deformation is not reported in their paper. According to \citet{Yuan2023}, a standard F-bar method is applied. From our investigation, the standard F-bar method is not able to efficiently overcome the stress oscillation in the case of a very large deformation {as shown in Fig.~\ref{2MP_MCslope_compareStab} (see also \citet{Xie2023b}). The standard F-bar method eliminates the volumetric locking under moderate deformation at $t = 2 \,s$. However, with the development of large deformation, a new sort of stress oscillation has been developed and enlarged as the propagation of the failure. In contrast, the proposed modified F-bar method shows promising performance in stabilising the suspicious stress oscillation at all scales of deformation.} Compared with a similar slope problem conducted by \citet{Kularathna2021}, our modified F-bar approach also provides smoother results. 

As we can see from Fig. \ref{MCslope_result}, the accumulated plastic strain contours given by the low and high conductivity conditions are quite different. {Accumulated plastic strain is an internal variable of this Mohr-Coulomb model \citep{deSouza2008}, measuring the cumulative plasticity the material has undergone.} The corresponding contour can indicate the location of the shear band and the failure mechanism. As we can see in Fig. \ref{MCslope_result}, the slope with relatively low permeability tends to develop a shear band at a deeper level. Therefore, more material points are mobilised in this failure mechanism{, leading to a larger settlement. This is because the pore pressure generated during the failure dissipates slowly in this case, resulting in a lower soil's effective stress, as shown in Fig. \ref{MCslope_result}}. For the high conductivity case, the location of the shear band is shallower. However, a secondary shear band is developing at a deeper level, as shown in Fig. \ref{MCslope_result}, although this band is less critical than the shallower one.

\subsection{Landslides with Nor-Sand model}

In this section, we investigate landslides using a more advanced soil model, Nor-Sand. The previous Mohr-Coulomb model can only describe the plastic behaviour of the soil without taking into account the densification or dilation effects (volumetric strains), which is another important property that can significantly influence its behaviour and variations in the pore water pressures. The Nor-Sand model is based on the critical state framework, and the soil density can be modelled using the specific volume (or void ratio). Therefore, we study landslides with two different soil densities by assigning different initial specific volumes. We set the initial specific volume $v_0$ to 1.55 and 1.70 for the cases of dense and loose sand, respectively. Nor-Sand material parameters for Hostun sand, calibrated by \citet{Andrade2008}, are adopted in this research. These material parameters are summarised in Table \ref{tabNS}, and the calibration with the experiment is shown in Appendix \ref{secA1}.
\begin{table}[h]
	\caption{Parameters for Landslides with Nor-Sand model}\label{tabNS}%
	\begin{tabular}{@{}ll@{}}
		\toprule
		Parameter &Value  \\
		\midrule
		Water density $\rho_w [kg/m^3]$    & 1000 \\
		Soil grain density $\rho_s [kg/m^3]$    & 2650 \\
		Shear modulus $G[MPa]$    & 40 \\
		Swelling index $\kappa$    & 0.002 \\
            Reference specific volume $v_{c0}$    & 1.892 \\
            Compression index $\lambda$    & 0.02 \\
		Slope of the critical state line $M$    & 1.0 \\
		Yield function constant $N$    & 0.1 \\
            Plastic potential constant $\Bar{N}$    & 0.1 \\
		Hardening coefficient $h$    & 100 \\
            Multiplier for maximum plastic dilatancy $\alpha$    & -2.0 \\
            Initial specific volume $v_0$    & 1.55 \& 1.70 \\
		Initial porosity $n_0$    & 0.4 \\
		Initial hydraulic conductivity $k_0 [m/s]$    & $5\times10^{-3}$ \\
		Gravitational acceleration $g[m/s^2]$    & -9.81 \\
		\bottomrule
	\end{tabular}
\end{table}

\begin{figure}
	\begin{center}
		\includegraphics[width=0.48\textwidth]{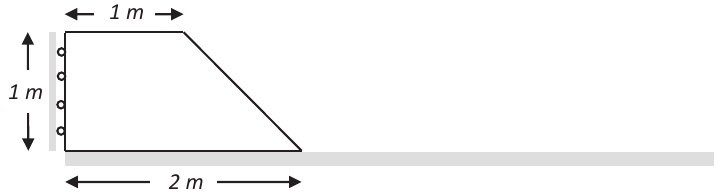}
		\caption{Landslides with Nor-Sand model: graphical illustration of the numerical model.}
		\label{NS}
	\end{center}
\end{figure}

Fig. \ref{NS} describes the geometry and boundary conditions for this numerical example. As in the previous cases, the roller and fixed boundary conditions are applied to the sides and bottom surfaces, respectively. {Similarly, impermeable boundaries are formed where those roller and fixed boundary conditions are applied, and the pore pressure is allowed to dissipate freely on the free water surface.} The critical state model has to be initialised to a site condition before proceeding with the analysis. The detailed initialisation procedure is presented in Appendix \ref{secA1}. During the initialisation, only the elastic model is used, and the gravity load is gradually applied in 1 s. Relatively high hydraulic conductivity (that is, $1\,m/s$) is assumed to better initialise pore pressures, and the water material points are attached to the soil material points to prevent seepage flow. {This procedure ensures the pore pressure generated during the initialisation is fully dissipated, representing an on-site pressure field before the landslides.} After the initialisation, the plasticity of the Nor-Sand model is suddenly enabled, and the landslide starts to be generated. The propagation of the landslides becomes nearly static at 12 s and 6 s (excluding the initialisation) for the cases of dense and loose sands, respectively. For this numerical example, a uniformly distributed $0.025\,m$ wide square grid with $3^2$ material points per cell is used, discretising the model with 21,540 material points for both the soil and the water phases. The results given by this configuration are in very good agreement with a configuration of refined mesh and an increased number of material points per cell {(i.e. $4^2$ material point per cell)}. The time step is set to $1.2\times10^{-4}\,s$, which is around 0.92 CFL. 

In this numerical example, the performance of single-point and double-point approaches is compared. From a computational cost point of view, the double-point approach is more cumbersome, but we find that this approach is only about 15\% slower than the single-point approach for the case of a quadratic B-spline shape function, which can be considered an acceptable increase of the computational effort. The additional computational effort is mainly coming from the construction of the shape function for the water material points. For the case of linear shape functions such as GIMPM, the differences in computational costs between the single and double-point approaches are even smaller because the construction of a linear shape function is very efficient.

\begin{figure}
	\begin{center}
		\includegraphics[width=0.8\textwidth]{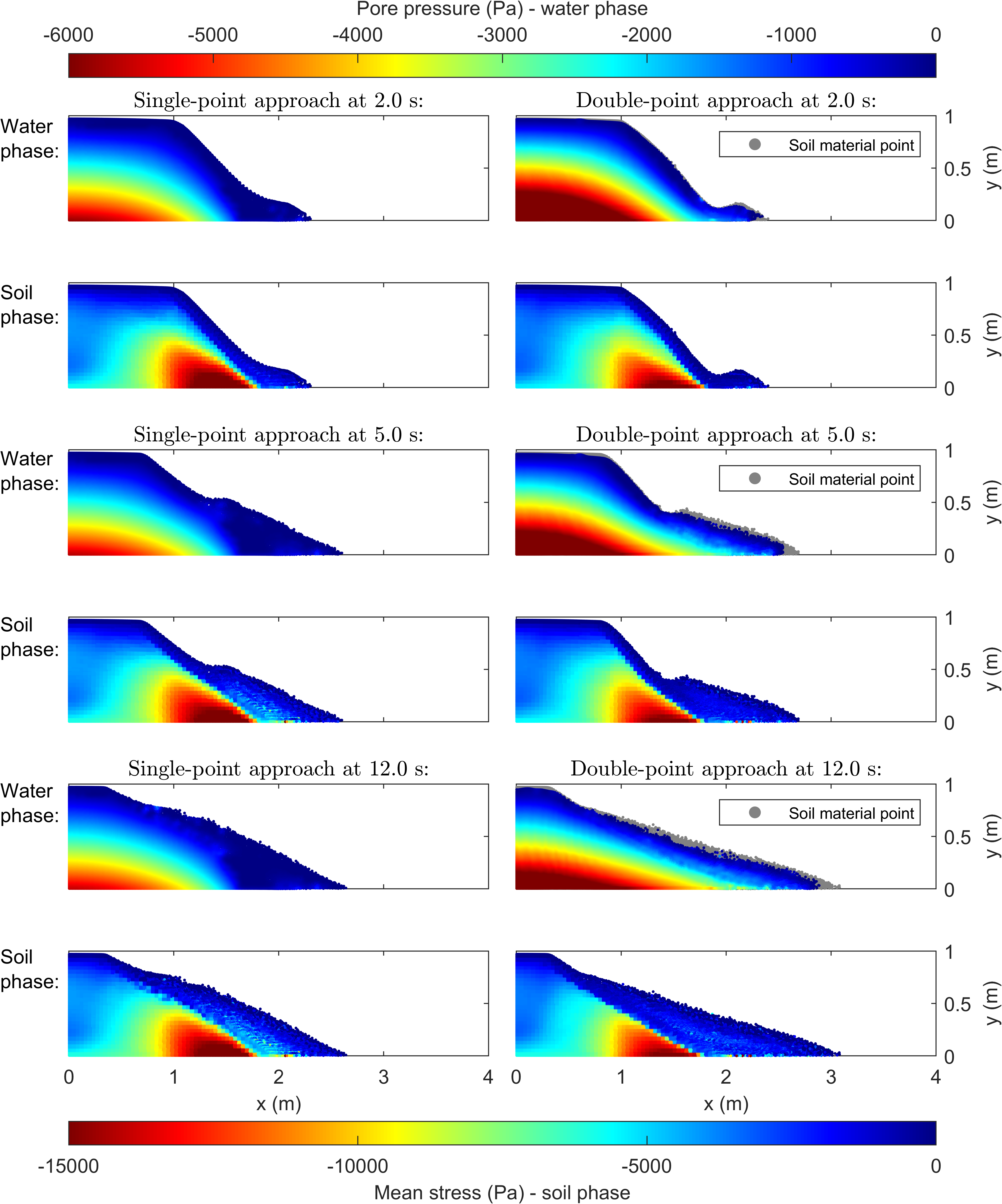}
		\caption{The evolution of landslides with Nor-Sand model (dense sand: $v_0=1.55$) - contours of pore water pressures and soil mean effective stresses given by single-point and double-point approaches}
		\label{NSdense}
	\end{center}
\end{figure}

Apart from computational efficiency, the numerical results obtained from the single-point and double-point approaches can sometimes be significantly different. As we can see in Fig. \ref{NSdense}, in the case of dense sand, the landslide simulations are significantly different between the single-point and double-point methods, and these differences increase as the deformations become greater. Landslides cannot be fully generated by the single-point method, resulting in a significantly underestimated propagation distance. {This significant difference is because the soil behaviour is dilative under dense conditions. The effective stress builds up during the dilation. In the double-point approach, the relative movement between the soil and water phases is accurately modelled. In this case, the pore water flow tends to drag the soil's dilation, resulting in a softer behaviour than the single-point approach. Due to the dilation, the volume of the soil phase is increased, while the water phase is incompressible, resulting in a small gap between the soil and water surfaces in the double-point method as shown in Fig. \ref{NSdense}. However, this phenomenon is not able to be captured by the single-point approach.} 
 
\begin{figure}
	\begin{center}
		\includegraphics[width=0.8\textwidth]{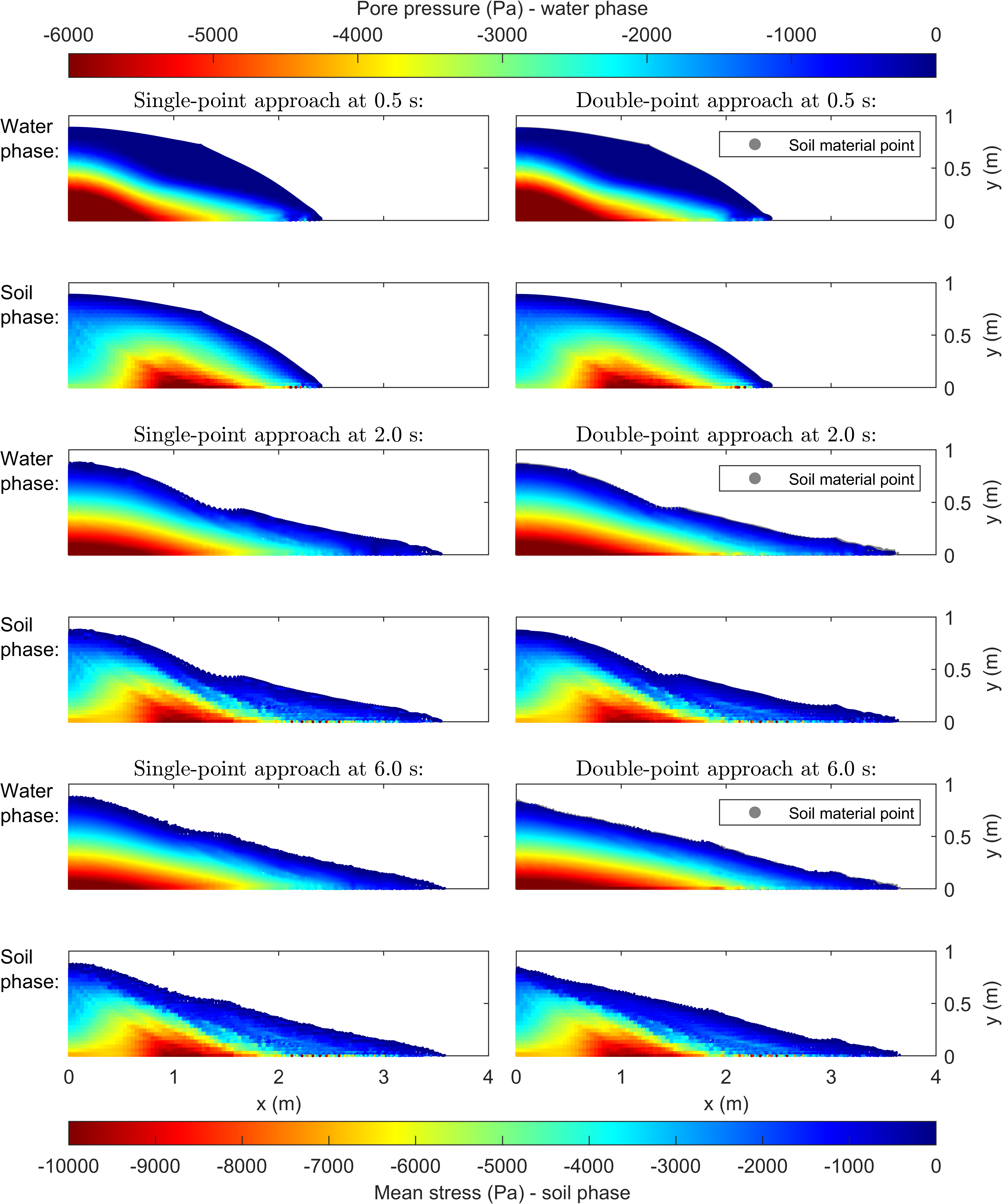}
		\caption{The evolution of landslides with Nor-Sand model (loose sand: $v_0=1.70$) - contours of pore water pressures and soil mean effective stresses given by single-point and double-point approaches}
		\label{NSloose}
	\end{center}
\end{figure}

Fig. \ref{NSloose} shows the comparison between single-point and double-point approaches under conditions of loose sand, for which the differences between the single-point and double-point approaches are smaller. In this case, the relative motion between the water and soil phases is not as high as in the case of dense sand. Instead, the water and soil material points are moving together rapidly in this failure mode. {The propagation distance is similar between the single-point approaches. However, the single-point approach still overestimates the stiffness on the slope crest after 2 seconds, as shown in Fig. \ref{NSloose}. In the loose condition, the soil is contractive in general. Therefore, only a very narrow gap is formed between the soil and water surfaces as shown in Fig. \ref{NSloose}. In the current formulation, the sand is assumed to be either fully dried or saturated. In other words, the proposed formulation cannot handle the suction that occurs in the partially saturated soil. However, in reality, suction exists in this gap between the soil and water surfaces depending on the degrees of saturation. The partially saturated formulation will be derived in future research. The suction effect may be insignificant in granular soil like sand because of its fast drainage.}

\begin{figure}
	\begin{center}
		\includegraphics[width=0.95\textwidth]{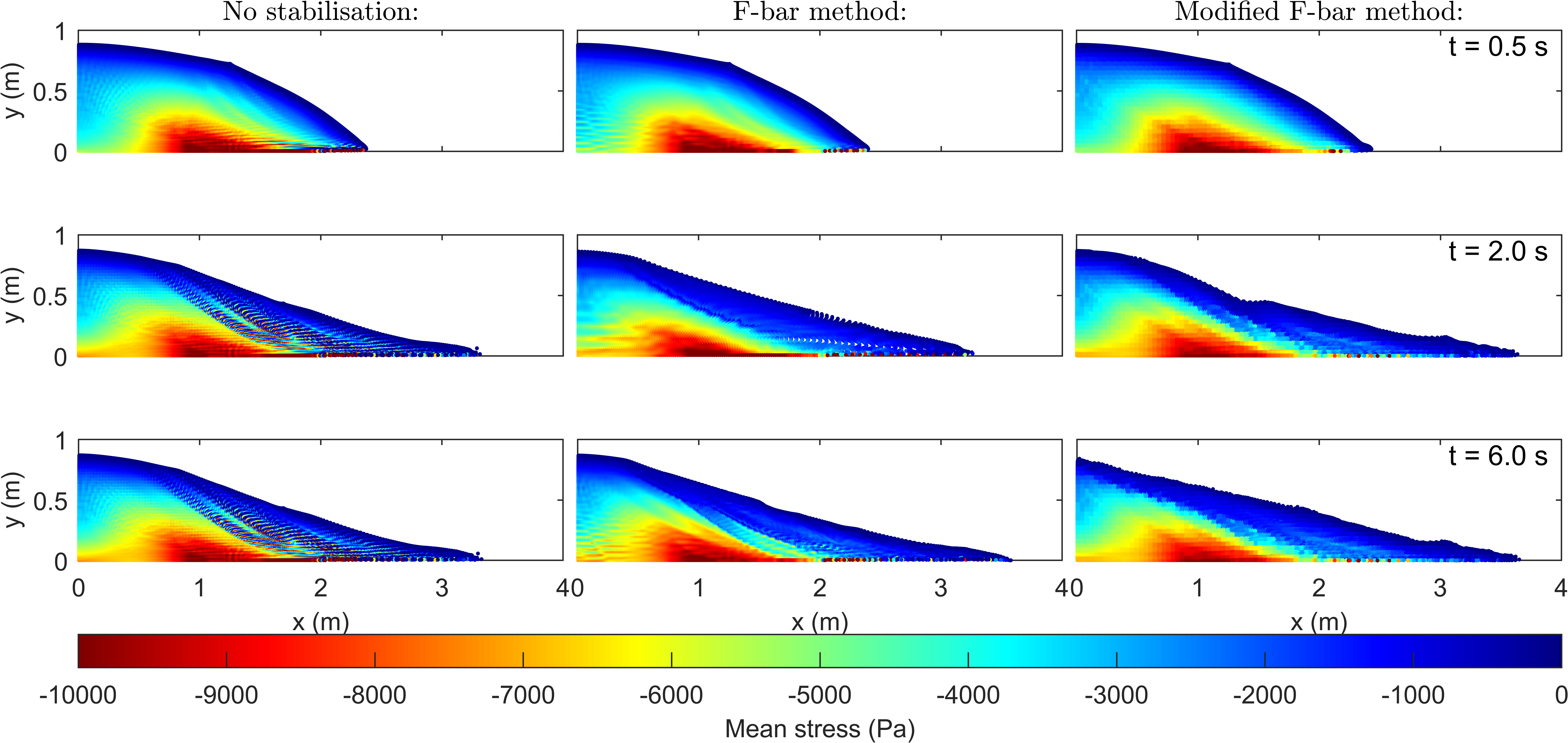}
		\caption{Soil mean effective stress contours of landslides with Nor-Sand model (loose sand: $v_0=1.70$) with different stabilisation methods at various time snapshots.}
		\label{2MP_NS_compareStab}
	\end{center}
\end{figure}

{As we can see in Figs. \ref{NSdense} and \ref{NSloose}, the pore pressure fields are very smooth and the stress fields are free of oscillation, indicating that the proposed stabilisation techniques are capable of efficiently stabilising the water and soil phases under very large deformations, such as the scenario of landslides. A very smooth stress field can be obtained by refining the mesh like the demonstration shown in Appendix \ref{secA:mFbar}. However, it is not necessary to this numerical problem according to the sensitivity analysis. The performance of the modified F-bar method has also been studied by comparing it with the case of no stabilisation and the F-bar method. As shown in Fig. \ref{2MP_NS_compareStab}, although the standard F-bar method has reduced the stress oscillation, there is still some oscillation at the bottom left of the slope. Also, an over-stiffening behaviour due to the volumetric locking can be observed in this numerical example. Compared to Fig. \ref{2MP_MCslope_compareStab}, the example using the Mohr-Coulomb model, the Nor-Sand constitutive model suffers less volumetric locking instability. This is because the flow rule of such an advanced constitutive model involves the evolution of the volumetric stain \citep{Borja2006}, resulting in less volumetric constraint. A similar phenomenon can also be observed in the research conducted by \citet{Xie2023} using the Modified Cam-Clay constitutive model. }

\section{Conclusion}\label{sec6}
In this research, a two-phase, double-point MPM has been derived based on the incremental fractional step method. Soil and water are considered solid and liquid for this research, respectively. The soil is assumed to be dry or fully saturated, and the water is assumed to be incompressible. For this coupled MPM, the pore pressure is solved implicitly. Therefore, the critical time step is governed by the CFL condition of the soil phase, and a relatively large time step can be used compared with a fully explicit approach. {Unlike the other formulations, the critical time step in the proposed formulation is not affected by the soil's permeability (or hydraulic conductivity), reducing the computational effort significantly in the case of low permeability.}

Cell-crossing noise is eliminated using B-spline shape functions. The mapping and remapping technique stabilises the oscillation in the water phase. Furthermore, we have improved the modified F-bar method to remove volumetric locking in the soil phase. Numerical examples show an excellent performance of the proposed stabilisation method. 

The proposed two-phase double-point MPM is validated with analytical solutions and numerical examples under small and large deformations, respectively. After that, the Nor-Sand constitutive model is used to study landslides. We highlight the importance of using a two-phase double-point MPM to study such a problem. {Single-point MPM cannot properly generate landslides when the soil is dilative (dense soil) because the relative displacement between the soil and water phase cannot be modelled properly. In this case, the separation between soil and water surface (a gap) formed by soil dilation cannot be captured by the single-point method. When the slope is initially loose, the single-point method also overestimates the stiffness around the crest.} Thus, the double-point approach is more accurate compared with the single-point method. Therefore, it is necessary to use a two-phase double-point MPM in the case of fully saturated soil to accurately simulate geotechnical engineering problems undergoing high deformations. Double-point MPM was criticised in the past as a cumbersome method. However, we demonstrate that it does not double the computational cost, and it only represents an increase of about 15\% of the calculation time compared with the single-point MPM if a quadratic B-spline shape function is used. This difference can even be reduced by improving the computational efficiency of constructing the shape function. 

{The limitation of the proposed method has been raised after the development of the gap (partially saturated zone) between the soil and water surfaces. The proposed method is not capable of handling the suction that exists in this zone. Therefore, it is necessary to further develop this approach to accommodate partially saturated soil. }

\begin{appendices}
\section{Elastic strip footing}\label{secA:mFbar}
{An elastic strip footing problem is used here to demonstrate the various performance of the cell-wise and patch-wise modified F-bar method on the GIMPM and BSMPM. The classic linear elastic constitutive model, the one used in the Mohr-Coulomb model, is used in this numerical example. The Young's modulus and Poisson's ratio are 10~$MPa$ and 0.49, respectively. In this case, the material is nearly incompressible. Roller boundary conditions are used for the right- and left-hand side boundaries to ensure a symmetric boundary condition. The elastic material has a geometry of 20 $m$ in width and 10 $m$ in height. Only half of the width is modelled due to the symmetric boundary condition. The plane strain condition is assumed for this problem. A 6 $m$ wide strip footing modelled by a 20 $kPa$ uniformly distributed load, spanning between 0 and 3 $m$ (6 $m$ wide in total because of the symmetric). The load is gradually applied within 1 $s$, and the gravitational acceleration is ignored. For this problem, a 0.2 $m$ square mesh ($h = 0.2\,m$) is used for the case of coarse mesh, and the mesh size is reduced to 0.05 $m$ to study the mesh size effect on the resolution of the stress contour. $2^2$ material points are placed in each cell for both cases. }

\begin{figure}
	\begin{center}
		\includegraphics[width=0.9\textwidth]{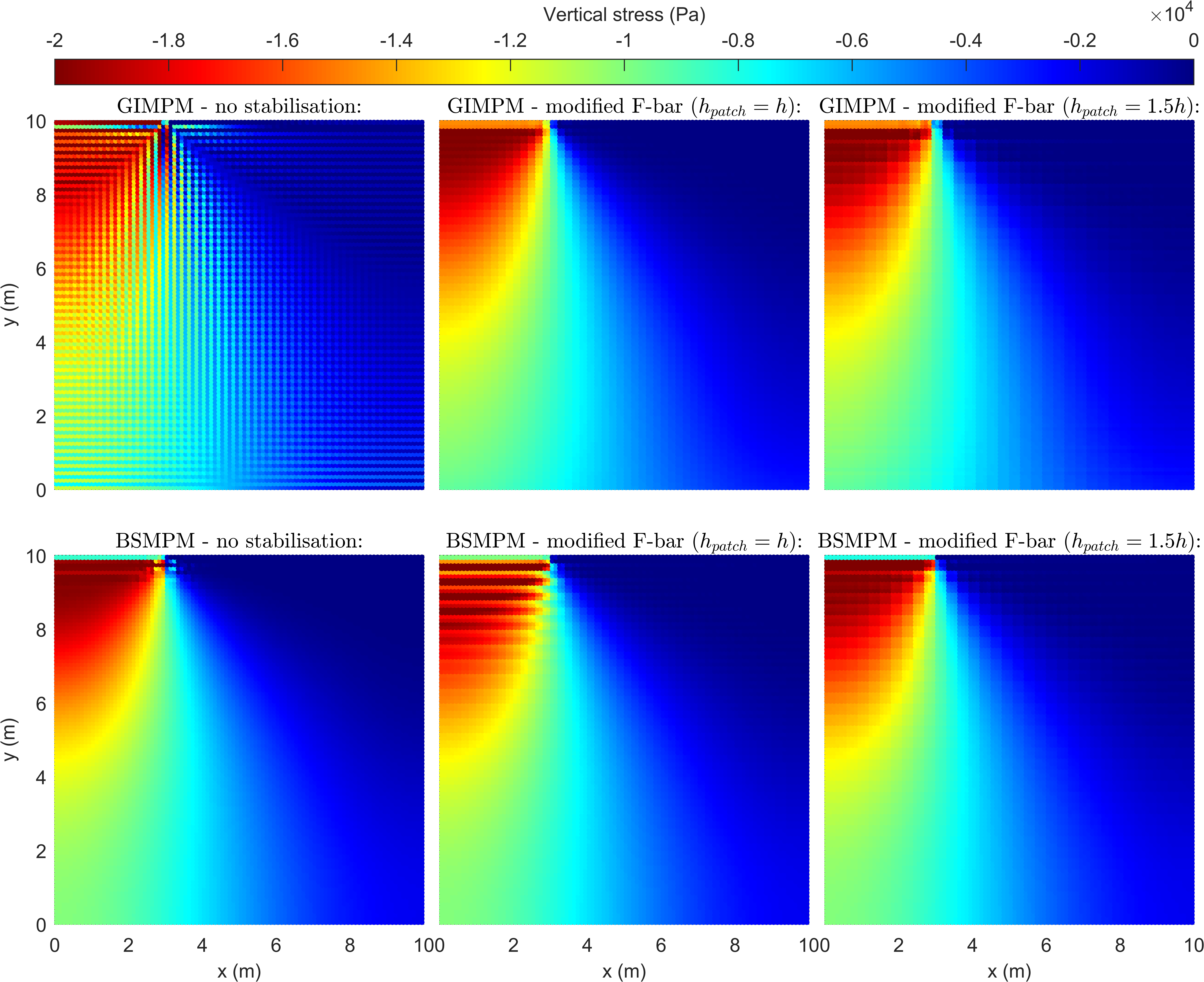}
		\caption{Comparison of the GIMPM and BSMPM with different stabilisation methods.}
		\label{Compare_stab_method_footing}
	\end{center}
\end{figure}

\begin{figure}
	\begin{center}
		\includegraphics[width=0.9\textwidth]{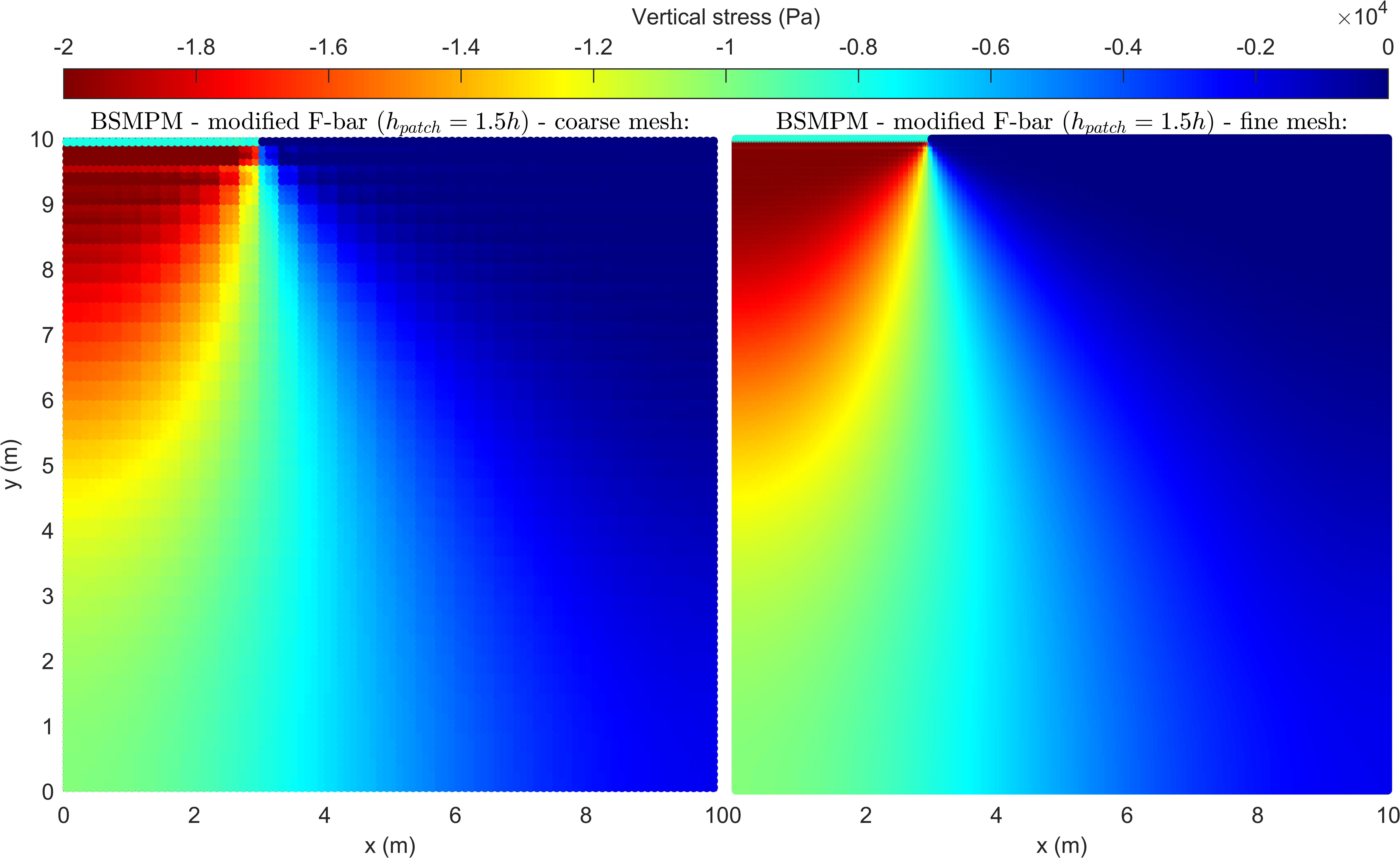}
		\caption{Illustration of the mesh size effect on the resolution for the stress contour given by BSMPM with modified F-bar method.}
		\label{Compare_mFbarPatch_meshSize}
	\end{center}
\end{figure}

{Fig.~\ref{Compare_stab_method_footing} compares the results given by GIMPM and BSMPM with no stabilisation, the modified F-bar with cell-wise ($h_{patch} = h$) and patch-wise ($h_{patch} = 1.5h$) averaging methods. We can see the GIMPM without stabilisation suffers strong volumetric locking. The higher-order B-spline shape function reduces the instabilities, but it cannot be fully removed without any stabilisation. However, applying the conventional cell-wise averaging ($h_{patch} = h$) on the BSMPM introduces new instability as shown in Fig.~\ref{Compare_stab_method_footing}, indicating this method is not valid for BSMPM. On the contrary, the modified F-bar method based on patch-wise averaging efficiently eliminates all the volumetric locking on both GIMPM and BSMPM. Finally, Fig.~\ref{Compare_mFbarPatch_meshSize} demonstrates the resolution of the stress contour can be improved by refining the mesh. }

\section{Nor-Sand model}\label{secA1}
\subsection{Numerical implementation}
\citet{Borja2006} show sufficient details of numerical implementation for their original implicit return-mapping algorithm. Therefore, only the modifications for the large deformation analysis have been provided in this paper. The notation in this section also follows \citet{Borja2006}.

\citet{Fern2016} reported that the implementation of the Nor-Sand model for a large deformation problem requires a dilatancy cutoff technique to prevent a sudden increase in dilatancy of a material point at a critical state. According to \citet{Fern2016}, the material points located on the free surface will be ejected from the flow body if the dilatancy suddenly increases during the shearing. For the implicit Nor-Sand algorithm, this instability will result in a non-physical complex number that diverges the Newton-Raphson iteration. Therefore, we introduce the same dilatancy cutoff technique as \citet{Fern2016} to improve the convergence performance for this implicit Nor-Sand algorithm. In this technique, a strain history variable, the accumulated deviatoric plastic strain $\Bar{E}^p_d$, is introduced, and the state parameter $\psi_i$ is calculated by
\begin{equation}\label{eq:eta}
	\begin{cases}
		\text{if}\; \Bar{E}^p_d \leq \Bar{E}^{p,\,start}_d : \\
		\psi_i = v - v_{c0} + \tilde{\lambda}\text{log}(-p_i) \\
		\\
		\text{if}\; \Bar{E}^{p,\,start}_d < \Bar{E}^p_d < \Bar{E}^{p,\,end}_d : \\
		\psi_i = (1-\frac{\Bar{E}^p_d - \Bar{E}^{p,\,start}_d}{\Bar{E}^{p,\,end}_d - \Bar{E}^{p,\,start}_d}) \psi_{(\Bar{E}^p_d = \Bar{E}^{p,\,start}_d)} \\
		\\
		\text{if}\; \Bar{E}^p_d \geq \Bar{E}^{p,\,end}_d : \\
		\psi_i = 0
	\end{cases} 
\end{equation}
where, $\Bar{E}^{p,\,start}_d$ and $\Bar{E}^{p,\,end}_d$ are user-defined variables. When $\Bar{E}^p_d \leq \Bar{E}^{p,\,start}_d$ the state parameter $\psi_i$ is updated as usual. After $\Bar{E}^p_d > \Bar{E}^{p,\,start}_d$ the residual state parameter will gradually reduce to zero and drop to nil at $\Bar{E}^p_d = \Bar{E}^{p,\,end}_d$. As in \citet{Fern2016}, we set $\Bar{E}^{p,\,start}_d = 100\%$ and $\Bar{E}^{p,\,end}_d = 150\%$ to ensure that it only applies to a sheared soil at a critical state.

As we know, the dry granular flow may experience a phase transformation from solid to gas, depending on its propagation speed. The critical state sand model can capture the phase transformation from solid (before the critical state) to liquid (reached the critical state). However, it is not intended to model a gas-like flow. In the $\mu$(I) rheology model \citep{Dunatunga2015}, a gas-like behaviour is achieved by setting the stress to zero when the material's density is lower than a critical value. In this case, the sand particles are disconnected. In the Nor-Sand model, the specific volume (or void ratio) represents the density of the sand. Therefore, the gas-like behaviour can be modelled by setting the stress to zero when the specific volume is larger than a critical value:
\begin{equation}
	\pmb{\sigma} = \pmb{0} \;\;\;\; \text{if} {(v>v_{crit})}.
\end{equation}
The corresponding elastic strain should also be nil. However, the left Cauchy green nil elastic tensor $\pmb{B}^e$ cannot be an exact identity matrix because an eigendecomposition must be performed on this tensor. Therefore, small, unequal values are assigned to its diagonal components if $v>v_{crit}$. As a result, the elastic logarithmic strain is approximately zero
\begin{equation}
	\pmb{E}^e = \frac{1}{2}\ln{\pmb{B}^e}\approx \pmb{0}.
\end{equation}
The critical specific volume $v_{crit} = 1.6v_{0}$ is used in this research. This 1.6 factor is consistent with \citet{Dunatunga2015}. A specific volume greater than the maximum specific volume means that the sand particles are disconnected, resulting in a nil stress and gas-like behaviour. The phase transformation from solid to gas is reversible. Once the specific volume falls below the defined critical value, the stress will be updated again by the Nor-Sand constitutive model. It is worth mentioning that, in the case of saturated soil, the mixture will behave like water instead of gas under the nil stress condition.

In the Nor-Sand algorithm \citet{Borja2006} derived, the specific volume of soil is updated by
\begin{equation}
    v^{t+1} = {J}^{t+1}\,v_{0} \,.
\end{equation}
However, we find that the convergence of this constitutive model under large deformation can be significantly improved by adopting
\begin{equation}
    v^{t+1} = \Bar{J}^{t+1}\,v_{0} \,,
\end{equation}
where
\begin{equation}
    \Bar{J}^{t+1} = \sum_{I=1}^{{N}_{n}} {N}_{Isp} \left(\frac{1}{m_{sI}}\sum_{p=1}^{{N}_{sp}} {N}_{Isp} \, {m}_{sp} \, {J}^{t+1}\right)
\end{equation}
is the stabilised Jacobian of deformation gradient. The difference in results is negligible between these two update procedures, and the performance of the constitutive model under large deformation has been significantly improved. 

\subsection{Single-element driver}

The single-element (or material point) test can be used to validate the implementation of the constitutive model. It is also an important tool for calibrating the model parameters with a stress-strain curve. The critical state constitutive model is usually calibrated with an undrained or drained triaxial test. Therefore, it is necessary to establish a single-element driver to simulate triaxial shear under undrained and drained conditions. In this section, we will derive the single-element driver for the constitutive models based on both infinitesimal and finite strain theories. 

\subsubsection{Undrained condition}
It is relatively easy to derive an undrained triaxial single-element driver. Under the undrained condition, the volume change is prevented, resulting in a stationary Jacobian $J=1$. In the case of an infinitesimal stain, the constitutive model is driven by the incremental strain. Therefore, we can establish a constant infinitesimal strain tensor $\pmb{\epsilon}$ by preventing its volume change
\begin{equation}\label{eq:strain_undrain}
	\pmb{\epsilon} = 
	\begin{bmatrix}
		\epsilon_a & 0 & 0 \\
		0 & -\epsilon_a/2 & 0 \\
		0 & 0 & -\epsilon_a/2 
	\end{bmatrix},
\end{equation}
where, $\epsilon_a$ is the infinitesimal axial strain, and the minus sign means compression. In this case, the Jacobian $ J=1+\text{tr}\,\pmb{\epsilon}=1$. Therefore, the incremental form of (\ref{eq:strain_undrain}) which drives the constitutive model can be written as
\begin{equation}
	\Delta\pmb{\epsilon} = 
	\begin{bmatrix}
		\Delta\epsilon_a & 0 & 0 \\
		0 & -\Delta\epsilon_a/2 & 0 \\
		0 & 0 & -\Delta\epsilon_a/2 
	\end{bmatrix}.
\end{equation}
In the case of finite strain, the axial stretch
\begin{equation}\label{eq:axial_stretch}
	\lambda_a = \frac{l}{l_0}
\end{equation}
is defined together with the logarithmic axial strain
\begin{equation}\label{eq:log_strain_a}
	\pmb{E}_a = \text{ln}\,\lambda_a = \text{ln}(l/l_0) \;,
\end{equation}
where, $l_0$ and $l$ are the length of the bar at reference and deformed configurations, respectively. We can replace the infinitesimal axial strain in (\ref{eq:strain_undrain}) with the logarithmic axial strain (\ref{eq:log_strain_a}) because the split of the logarithmic axial strain is also additive. Therefore, Equation (\ref{eq:strain_undrain}) yield to
\begin{equation}
	\pmb{E} = 
	\begin{bmatrix}
		\text{ln}(l/l_0) & 0 & 0 \\
		0 & \text{ln}(l_0/l)/2 & 0 \\
		0 & 0 & \text{ln}(l_0/l)/2 
	\end{bmatrix},
\end{equation}
We know the relationship between the left Cauchy-Green tensor $\pmb{B}$ and logarithmic strain 
\begin{equation}
	\pmb{E} = \frac{1}{2} \text{ln}\,\pmb{B} \;.
\end{equation}
Then, we can easily find the left Cauchy-Green tensor under the undrained triaxial condition by 
\begin{equation}
	\pmb{B} = \text{exp}(2\pmb{E}) =
	\begin{bmatrix}
		{(l/l_0)}^2 & 0 & 0 \\
		0 & l_0/l & 0 \\
		0 & 0 & l_0/l 
	\end{bmatrix} .
\end{equation}
The deformation gradient for this undrained triaxial extension or compression can be obtained by
\begin{equation}
	\pmb{F} = \sqrt{\pmb{B}} = \text{exp}\,\pmb{E} = 
	\begin{bmatrix}
		l/l_0 & 0 & 0 \\
		0 & \sqrt{l_0/l} & 0 \\
		0 & 0 & \sqrt{l_0/l} 
	\end{bmatrix} .
\end{equation}
In the case of finite strain, the incremental deformation gradient $\Delta \pmb{F}$ drives the constitutive model. Therefore, by assuming $l_0=1$ and $l = l_0 + \Delta\epsilon_a$ we can have
\begin{equation}\label{eq:dF_undrained}
	\Delta\pmb{F} =
	\begin{bmatrix}
		1 + \Delta\epsilon_a & 0 & 0 \\
		0 & \sqrt{1/(1 + \Delta\epsilon_a)} & 0 \\
		0 & 0 & \sqrt{1/(1 + \Delta\epsilon_a)} 
	\end{bmatrix} .
\end{equation}
Also, $J = \text{det}\,\pmb{F} = \text{det}\,\Delta\pmb{F} = 1$ because of the prevented volume change under the undrained condition. 
\subsubsection{Drained condition}
It is also possible to derive a single-element drained triaxial driver. We know the stress path in a drained triaxial test always follows
\begin{equation}
	\Delta q = 3\Delta p \,.
\end{equation}
To maintain this condition under both elastic and plastic loading, we have to relate the incremental strain with the elastoplastic tangent modulus $\mathbb{C}^{ep}$. Conveniently, the Voigt notation for this fourth-order tensor ${\text{C}}_{ij}$ is used. Therefore, the incremental infinitesimal strain for the drained triaxial driver can be written as
\begin{equation}
	\Delta\pmb{\epsilon} = 
	\begin{bmatrix}
		\Delta\epsilon_a & 0 & 0 \\
		0 & -\Delta\epsilon_a \frac{\text{C}_{21}}{\text{C}_{22}+\text{C}_{23}}
		 & 0 \\
		0 & 0 & -\Delta\epsilon_a\frac{\text{C}_{31}}{\text{C}_{32}+\text{C}_{33}}
	\end{bmatrix}.
\end{equation}
Following the same procedures from (\ref{eq:axial_stretch}) to (\ref{eq:dF_undrained}), we can derive the deformation gradient for the drained triaxial driver as
\begin{equation}
	\pmb{F} =   
	\begin{bmatrix}
		l/l_0 & 0 & 0 \\
		0 & ({l_0/l})^{\frac{\text{C}_{21}}{\text{C}_{22}+\text{C}_{23}}} & 0 \\
		0 & 0 & ({l_0/l})^{\frac{\text{C}_{31}}{\text{C}_{32}+\text{C}_{33}}} 
	\end{bmatrix} ,
\end{equation}
and its incremental form
\begin{align}
	\Delta\pmb{F} =   
	\begin{bmatrix}
		1 + \Delta\epsilon_a & 0 & 0 \\
		0 & (\frac{1}{1 + \Delta\epsilon_a})^{\frac{\text{C}_{21}}{\text{C}_{22}+\text{C}_{23}}} & 0 \\
		0 & 0 &  (\frac{1}{1 + \Delta\epsilon_a})^{\frac{\text{C}_{31}}{\text{C}_{32}+\text{C}_{33}}}
	\end{bmatrix} .
\end{align}

\subsection{Calibration and validation}

%%%%%%%%%%
\begin{figure*}
	\begin{center}
		\includegraphics[width=\textwidth]{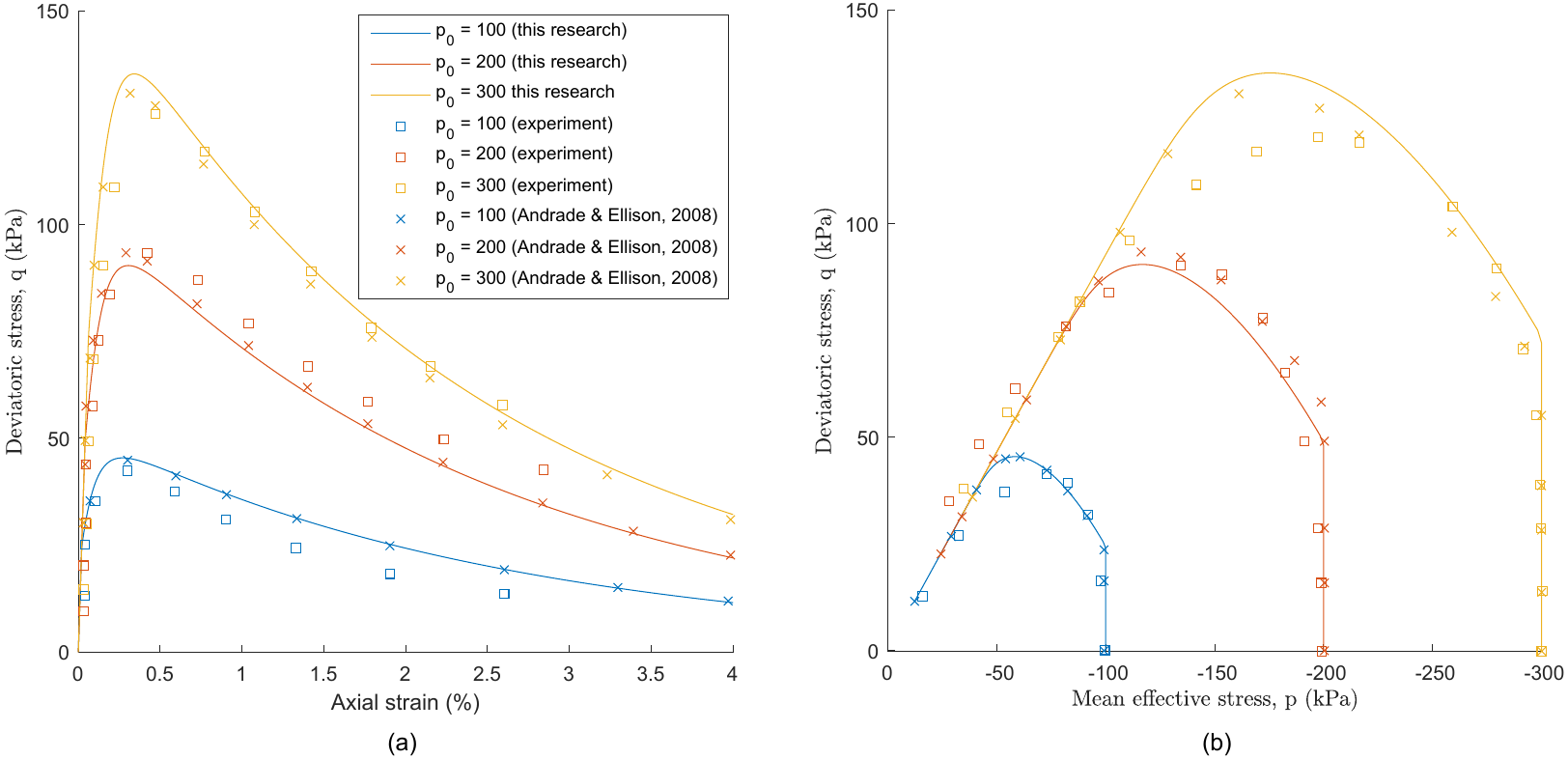}
		\caption{Calibration and validation of Nor-Sand constitutive model versus undrained triaxial experimental data and the corresponding numerical study}
		\label{undrained_validation}
	\end{center}
\end{figure*}
%%%%%%%%%%

\citet{Andrade2008} calibrated the Nor-Sand model parameters for Hostun sand based on the undrained triaxial experiments conducted by \citet{Doanh1997}. Table \ref{tabHostunSand} summarises the calibrated parameters for Hostun sand. Our numerical implementation for the Nor-Sand model uses the single-element driver derived in the previous section. Three initial mean effective stresses $p_0$, ranging from 100 to 300 kPa, are used. The initial specific volume is 2.0, and the initial elastic volumetric strain is 0. Fig. \ref{undrained_validation} shows the validation results. As we can see, the results of our implementation are closely correlated with \citet{Andrade2008} and also the experiments conducted by \citet{Doanh1997}, which implies the correctness of the implementation of the constitutive model.   

%%%%%%%%%%
\begin{table}[h]
	\caption{Calibrated parameters for Hostun sand specimens}\label{tabHostunSand}%
	\begin{tabular}{@{}ll@{}}
		\toprule
		Parameter &Value  \\
		\midrule
		Shear modulus $G[MPa]$    & 40 \\
		Swelling index $\kappa$    & 0.002 \\
            Reference specific volume $v_{c0}$    & 1.892 \\
            Compression index $\lambda$    & 0.02 \\
		Slope of the critical state line $M$    & 1.0 \\
		Yield function constant $N$    & 0.1 \\
            Plastic potential constant $\Bar{N}$    & 0.1 \\
		Hardening coefficient $h$    & 100 \\
            Multiplier for maximum plastic dilatancy $\alpha$    & -2.0 \\
		\bottomrule
	\end{tabular}
\end{table}
%%%%%%%%%%

\subsection{Initialise the multi-element MPM model}
In the Mohr-Coulomb/Drucker-Prager model, the soil's strength is related to the friction angle and cohesion. Also, constant Young's modulus and Poisson's ratio are used for the elastic behaviour. In this case, its strength does not depend on the initial condition. For this hyperelastic Nor-Sand model, apart from the parameters that can be calibrated from the triaxial test, we need to initialise the reference pressure $p_0$, the elastic volumetric strain at the reference pressure $\epsilon^e_{v0}$, and the initial image pressure.

It is logical to obtain the initial condition by gradually imposing the gravity load. To prevent failure during gravity loading, only elastic behaviour is allowed at this stage. Initially, we set the reference pressure $p_0 = 2$ kPa and the elastic volumetric strain $\epsilon^e_{v0} = 0$ at this reference pressure. Then, the gravity load is gradually applied from zero to its maximum value. During gravity loading, the elastic stress and strain are updated by the hyperelastic law. Therefore, we can obtain the pressure and elastic volumetric strain by
\begin{equation}
	p = \frac{1}{3}\text{tr}\,\pmb{\sigma} \,,
\end{equation}
and
\begin{equation}
	\epsilon^e_v = \text{tr}\,\pmb{\epsilon}^e \,.
\end{equation}
At the end of gravity loading, we can take the last converged $p$ and $ \epsilon^e_v$ as reference pressure $p_0$ and elastic volumetric strain $\epsilon^e_{v0}$. From the yield function of the Nor-Sand model, we can obtain the initial image pressure $p_i$ using an over-consolidation ratio (OCR) and the reference pressure $p_0$ by
\begin{equation}\label{eq:pi_p0}
	p_i = OCR \times
	\begin{cases}
		p_0/\text{exp}(1) \\
		p_0(1-N)^{(1-N)/N}
	\end{cases} \begin{array}{c}
		\text{if}\; N=0, \\
		\text{if}\; N>0.
	\end{array}
\end{equation}
Also, the initial specific volume $v_0$ after the elastic loading can be updated by
\begin{equation}
	v_0 \leftarrow Jv_0 \,.
\end{equation}
Note that the volume change $J$ is very small during elastic loading, resulting in a minor change in the specific volume. However, imposing a well-gradient initial specific volume (i.e. a higher specific volume at the upper layer) can help the algorithm's convergence.

%An appendix contains supplementary information that is not an essential part of the text itself but which may be helpful in providing a more comprehensive understanding of the research problem or it is information that is too cumbersome to be included in the body of the paper.




\end{appendices}

\bibliography{sn-bibliography}% common bib file

\end{document}